\documentclass[10pt]{amsart}
\usepackage{amssymb}

\theoremstyle{definition}
\newtheorem{thm}{Theorem}[section]
\newtheorem{lem}[thm]{Lemma}
\newtheorem{prp}[thm]{Proposition}
\newtheorem{dfn}[thm]{Definition}
\newtheorem{cor}[thm]{Corollary}
\newtheorem{cnj}[thm]{Conjecture}

\newtheorem{exa}[thm]{Example}

\newenvironment{pff}{{\emph{Proof:}}}{\QED}

\newcommand{\bit}{\begin{itemize}}
\newcommand{\eit}{\end{itemize}}
\newcommand{\beq}{\begin{equation}}
\newcommand{\eeq}{\end{equation}}
\newcommand{\beqr}{\begin{eqnarray*}}
\newcommand{\eeqr}{\end{eqnarray*}}
\newcommand{\bal}{\begin{align*}}
\newcommand{\eal}{\end{align*}}
\newcommand{\ts}{\textstyle}
\newcommand{\rsz}[1]{\raisebox{0ex}[0.8ex][0.8ex]{$#1$}}

\newcommand{\af}{\alpha}
\newcommand{\bt}{\beta}
\newcommand{\gm}{\gamma}
\newcommand{\dt}{\delta}
\newcommand{\ep}{\varepsilon}
\newcommand{\zt}{\zeta}
\newcommand{\et}{\eta}
\newcommand{\ch}{\chi}
\newcommand{\io}{\iota}
\newcommand{\te}{\theta}
\newcommand{\ld}{\lambda}
\newcommand{\sm}{\sigma}

\newcommand{\ph}{\varphi}
\newcommand{\ps}{\psi}
\newcommand{\rh}{\rho}
\newcommand{\om}{\omega}
\newcommand{\ta}{\tau}

\newcommand{\Q}{{\mathbf{Q}}}
\newcommand{\Z}{{\mathbf{Z}}}
\newcommand{\R}{{\mathbf{R}}}
\newcommand{\C}{{\mathbf{C}}}
\newcommand{\N}{{\mathbf{N}}}

\newcommand{\tsr}{{\mathrm{tsr}}}

\newcommand{\Tr}{{\mathrm{Tr}}}
\newcommand{\id}{{\mathrm{id}}}
\newcommand{\ev}{{\mathrm{ev}}}
\newcommand{\sint}{{\mathrm{int}}}

\newcommand{\rank}{{\mathrm{rank}}}

\newcommand{\sa}{{\mathrm{sa}}}
\newcommand{\spec}{{\mathrm{sp}}}

\newcommand{\dirlim}{\displaystyle \lim_{\longrightarrow}}

\newcommand{\Mi}{M_{\infty}}

\newcommand{\andeqn}{\,\,\,\,\,\, {\text{and}} \,\,\,\,\,\,}
\newcommand{\QED}{\rule{0.4em}{2ex}}

\newcommand{\ca}{C*-algebra}
\newcommand{\ct}{continuous}
\newcommand{\pj}{projection}
\newcommand{\mf}{manifold}

\newcommand{\diff}{diffeomorphism}

\newcommand{\nbhd}{neighborhood}
\newcommand{\cpt}{compact Hausdorff}
\newcommand{\hm}{homomorphism}
\newcommand{\wolog}{without loss of generality}
\newcommand{\Wolog}{Without loss of generality}

\newcommand{\ifo}{if and only if}
\newcommand{\rsha}{recursive subhomogeneous algebra}
\newcommand{\Rsha}{Recursive subhomogeneous algebra}

\newcommand{\lscp}{lower semicontinuous projection}
\newcommand{\rshd}{recursive subhomogeneous decomposition}

\newcommand{\tdim}{topological dimension}

\newcommand{\sdg}{slow dimension growth}
\newcommand{\ndg}{no dimension growth}
\newcommand{\mh}{minimal homeomorphism}

\title[Limits of recursive subhomogeneous algebras]{Cancellation and
stable rank for direct limits of recursive subhomogeneous
algebras}

\author{N.\  Christopher Phillips}

\address{Department of Mathematics, University  of Oregon,
       Eugene OR 97403-1222, USA}

\thanks{Research
      partially supported by NSF grants DMS 9400904 and DMS 9706850.}

\begin{document}

\setcounter{section}{-1}

\begin{abstract}

We prove the following results for a unital simple direct limit $A$ of
\rsha s with no dimension growth:

(1)
$\tsr (A) = 1$.

(2)
The \pj s in $\Mi (A)$ satisfy cancellation:
if $e \oplus q \sim f \oplus q$, then $e \sim f$.

(3)
$A$ satisfies
Blackadar's Second Fundamental Comparability Question:
if $p, \, q \in \Mi (A)$ are \pj s such that $\ta (p) < \ta (q)$
for all normalized traces $\ta$ on $A$, then $p \precsim q$.

(4)
$K_0 (A)$ is unperforated for the strict order:
if $\et \in K_0 (A)$ and there is $n > 0$ such that
$n \et > 0$, then $\et > 0$.

The last three of these results hold under certain weaker dimension
growth conditions and without assuming simplicity.
We use these results to obtain previously unknown information on the
ordered K-theory of the crossed product $C^* (\Z, X, h)$ obtained
from a \mh\  of an infinite finite dimensional compact metric space $X$.
Specifically, $K_0 (C^* (\Z, X, h))$ is unperforated for the strict
order,
and satisfies the following K-theoretic version of
Blackadar's Second Fundamental Comparability Question:
if $\et \in K_0 (A)$ satisfies $\ta_* (\et) > 0$
for all normalized traces $\ta$ on $A$, then there is a
\pj\   $p \in \Mi (A)$ such that $\et = [p]$.

\end{abstract}

\maketitle

\section{Introduction}

\Rsha s were introduced in \cite{PhN}; we recall the definition below.
They include finite direct sums of not necessarily trivial
unital homogeneous \ca s,
the dimension drop intervals and matrix algebras over them,
and the algebras $A_Y$ arising in Qing Lin's study \cite{Ln} of the
transformation group \ca s of minimal homeomorphisms (provided
$\sint (Y) \neq \varnothing$).
In this paper, we generalize to certain direct limits of
\rsha s some of the known results on direct limits of
homogeneous \ca s with slow dimension growth.
We use these results to obtain previously unknown information
on the ordered K-theory of the \ca s of \mh s.
In particular, our results make it possible in many cases to compute the
Elliott invariant \cite{El5} for the crossed product by a \mh.

We prove the following results for direct limits of \rsha s.
(See later in the introduction, and Section~1, for explanations
of the dimension growth conditions.)

\begin{thm}\label{C01}
Let $A = \dirlim (A_i, \ph_{i j})$ be a unital direct limit
of a system of \rsha s with slow dimension growth.
Then:
\bit
\item[(1)]
The map $U (A) / U_0 (A) \to K_1 (A)$ is an isomorphism.
\item[(2)]
If in addition the system has strict slow dimension growth,
then the \pj s in $\Mi (A)$ satisfy cancellation:
if $e \oplus q \sim f \oplus q$, then $e \sim f$.
\item[(3)]
If in addition the maps $\ph_{i j}$
of the system are all injective and unital,
then $A$ satisfies
Blackadar's Second Fundamental Comparability Question
(\cite{Bl3}, 1.3.1):
if $p, \, q \in \Mi (A)$ are \pj s such that $\ta (p) < \ta (q)$
for all normalized traces $\ta$ on $A$, then $p \precsim q$.
\item[(4)]
If in addition the maps $\ph_{i j}$
of the system are all injective and unital, and $A$ is simple,
then $K_0 (A)$ is unperforated for the strict order.
That is, if $\et \in K_0 (A)$ and there is $n > 0$ such that
$n \et > 0$, then $\et > 0$.
\item[(5)]
If in addition the system has no dimension growth,
and $A$ is simple, then $\tsr (A) = 1$.
\eit
\end{thm}

In this theorem, and throughout the paper, notation is as follows.
The notation $p \sim q$ means Murray-von Neumann equivalence of \pj s,
$p \precsim q$ means $p$ is Murray-von Neumann equivalent to a
sub\pj\  of $q$,
$\Mi (A)$ is the {\emph{algebraic}} direct limit
$\dirlim M_n (A)$ under the maps $a \mapsto a \oplus 0$
(following Definition~5.1.1 of \cite{Bl2}),
$\tsr (A)$ is the topological stable rank of $A$ (\cite{Rf1}),
$U (A)$ is the unitary group of a unital \ca\  $A$, and
$U_0 (A)$ is the identity component of $U (A)$.

Parts~(2), (3), and~(4) of the theorem partially generalize results of
\cite{Bl4}, \cite{MP}, and \cite{Gd}.
Part~(5) generalizes \cite{DNNP}.
We have had to impose extra conditions in part~(2)
(strict \sdg) and part~(4) (simplicity); we do not know whether
these extra conditions are really necessary.
Similarly, we do not know whether the condition
``\ndg'' in Part~(5) can be relaxed to ``\sdg'', as was
done for the homogeneous case in \cite{BDR}.

Two other results from the homogeneous case,
Theorem~2 of \cite{BDR} and Theorem~2.7 of \cite{Gd},
do not generalize to direct limits of \rsha s.
Specifically, there is a simple
direct limit $A$ of separable \rsha s, with \ndg\  and in which
the maps of the system are all injective and unital,
such that the projections in $A$ distinguish the traces on $A$
but $A$ does not have real rank zero, and such that
$K_0 (A)$ does not have Riesz decomposition.
Examples of these sorts were known before, but we give
one (Example~\ref{Q9}) which arises naturally from our
applications to crossed products by \mh s and has a simple proof.
(In \cite{Ph8}, we will give criteria for exactly when a
simple direct limit of separable \rsha s with \ndg\  has real
rank zero, and when it has the property (SP): every nonzero
hereditary subalgebra contains a nonzero \pj.
In particular, we will see that the combination of
Riesz decomposition, (SP), and projections distinguish traces,
implies real rank zero.)

Let $X$ be an infinite compact metric space,
and let $h$ be a minimal homeomorphism of $X$.
The crossed product $C^* (\Z, X, h)$ has been well studied
when $X$ is the Cantor set (see for example \cite{GPS}),
and when $h$ is an irrational rotation of the circle.
The \ca s of \mh s of higher dimensional spaces have remained
somewhat mysterious.
Connes has shown \cite{Cn} that, unlike the cases
mentioned above, the crossed products sometimes have
no nontrivial \pj s.
(See Corollary~12 in Section~6 of \cite{Ex} for a generalization.)
Qing Lin \cite{Ln} has studied simple subalgebras of the crossed product
which can be realized as direct limits of \rsha s in which
the maps of the system are all injective and unital;
moreover, if $X$ is finite dimensional, then
the system has no dimension growth.
Using those subalgebras and the results above, we obtain
the following theorem.
(Here, $U (A)$ is the unitary group of a \ca\  $A$, and
$U_0 (A)$ is the identity component of $U (A)$.)

\begin{thm}\label{Q05}
Let $X$ be a finite dimensional infinite compact metric space,
and let $h$ be a minimal homeomorphism of $X$.
Then:
\bit
\item[(1)]
The map
\[
U (C^* (\Z, X, h)) / U_0 (C^* (\Z, X, h))
                 \longrightarrow K_1 (C^* (\Z, X, h))
\]
is surjective.
\item[(2)]
$C^* (\Z, X, h)$ satisfies the following K-theoretic version of
Blackadar's Second Fundamental Comparability Question:
if $\et \in K_0 (A)$ satisfies $\ta_* (\et) > 0$
for all normalized traces $\ta$ on $A$, then there is a
\pj\   $p \in \Mi (A)$ such that $\et = [p]$.
\item[(3)]
$K_0 (C^* (\Z, X, h))$ is unperforated for the strict order.
\eit
\end{thm}

This theorem can be used to completely determine the order
on the $K_0$-group of the crossed product in interesting cases.
In Example~\ref{FT},
we easily obtain the description, proved in \cite{Kd},
of the positive cone in the $K_0$-group of the crossed product by a
Furstenberg transformation of the $2$-torus.
We also completely determine the Elliott invariant
(see \cite{El5})
for the crossed product by a \mh\  $h$ of an odd sphere $S^n$
with $n \geq 3$.
It follows from our computation that the Elliott invariant
depends only on the simplex of invariant Borel probability measures
for $h$, and in particular not on the dimension $n$ of the sphere
(as long as $n \geq 3$).
The Elliott classification conjecture would therefore imply that
if $n_1, \, n_2 \geq 3$ are odd, and $h_j$ is a uniquely ergodic
\mh\  of $S^{n_j}$,
then $C^* (\Z, S^{n_1}, h_1) \cong C^* (\Z, S^{n_2}, h_2)$.

We now recall the definition of a \rsha\  and some useful
associated terminology.
(See Definitions 1.1 and 1.2 of \cite{PhN}.)
First recall that if $A$, $B$, and $C$ are \ca s, and
$\ph \colon A \to C$ and $\rh \colon B \to C$ are \hm s, then
the pullback $A \oplus_C B$ is given by
\[
A \oplus_C B =
  \{ (a, b) \in A \oplus B \colon \ph (a) = \rh (b) \}.
\]

\begin{dfn}\label{A2}
A \rsha\  is a \ca\  of the form
\[
R = \left[ \cdots \rule{0em}{3ex} \left[ \left[
  C_0 \oplus_{C_1^{(0)}} C_1 \right]
 \oplus_{C_2^{(0)}} C_2 \right] \cdots \right]
            \oplus_{C_l^{(0)}} C_l,
\]
with $C_k = C \left( X_k, \, M_{n (k)} \right)$ for \cpt\  spaces $X_k$ and positive
integers $n (k)$, with
$C_k^{(0)} = C \left( \rsz{ X_k^{(0)}, \, M_{n (k)} } \right)$ for compact
subsets $X_k^{(0)} \subset X_k$ (possibly empty), and where the maps
$C_k \to C_k^{(0)}$ are always the restriction maps.
An expression of this type will be referred to as a
{\emph{decomposition}} of $R$ (over $\coprod_{k = 0}^l X_k$).

Associated with this decomposition are:
\bit
\item[(1)]
its {\emph{length}} $l;$
\item[(2)]
its {\emph{base spaces}} $X_0, X_1, \dots, X_l$ and
{\emph{total space}} $X = \coprod_{k = 0}^l X_k;$
\item[(3)]
its {\emph{matrix sizes}} $n (0), \dots, n (l)$, and
{\emph{matrix size function}} $m \colon X \to \N \cup \{0\}$,
defined by $m (x) = n (k)$ when $x \in X_k$ (this is called
the {\emph{matrix size of $A$ at $x$}});
\item[(4)]
its {\emph{minimum matrix size}} $\min_k n (k)$ and
{\emph{maximum matrix size}} $\max_k n (k);$
\item[(5)]
its {\emph{topological dimension}} $\dim (X)$
(the covering dimension of $X$ \cite{Pr}, Definition~3.1.1;
here equal to  $\max_k \dim (X_k)$), and
{\emph{topological dimension function}} $d \colon X \to \N \cup \{0\}$,
defined by $d (x) = \dim (X_k)$ when $x \in X_k$ (this is called
the {\emph{topological dimension of $A$ at $x$}});
\item[(6)]
its {\emph{standard representation}}
$\sm = \sm_R
    \colon R \to \bigoplus_{k = 0}^l C \left( X_k, \, M_{n (k)} \right)$,
defined by forgetting the restriction to a subalgebra in each of
the fibered products in the decomposition;
\item[(7)]
the associated
{\emph{evaluation maps}} $\ev_x \colon R \to M_{n (k)}$ for $x \in X_k$,
defined to be the restriction of the usual evaluation map to
$R$, identified with a subalgebra of
$\bigoplus_{k = 0}^l C \left( X_k, \, M_{n (k)} \right)$ via $\sm$.
\eit
\end{dfn}

At this point, we make a few remarks on the notions of \sdg.
For direct limits of direct sums of homogeneous \ca s,
it is usual to assume that the spaces associated with the
summands are connected.
Slow dimension growth is then defined in terms of the
dimensions of these spaces and the multiplicities of the
partial maps between direct summands at one level and
those at later levels.
(See  \cite{MP} and \cite{Gd}.)
Connectedness ensures that these multiplicities are well
defined.
In a direct system of \rsha s, because of the way the
summands ``overlap'' in the pullbacks, the multiplicities of the
partial maps between components of the \rshd s, and even the
partial maps themselves, need not be well defined.
Because of this, and because of what happens in some of our
proofs, it is not clear what the right definition of \sdg\  is.
A finite direct sum of algebras $C (X, M_n)$
is a single \rsha.
Therefore the right definition should include the situations
of \cite{MP} and \cite{Gd}.
It should also enable one to prove
cancellation in the general case, and
stable rank $1$ in the simple case.
Slow dimension growth might be more tractable for direct systems of
noncommutative CW-complexes \cite{Pd} and cell morphisms
(Definition~11.3 of \cite{Pd}).

We do not make a serious effort here to find the right
definition.
Rather, we give several versions which suffice for the
proofs of our theorems, and which apply to the algebras
we are most interested in, namely
simple direct limits with \ndg.
(These are the algebras required for the applications to
the \ca s of \mh s.)
We also don't formally consider weakenings of the dimension growth
conditions to ``relatively large entries'' in the sense of Section
3 of \cite{MP}.
(One can, however, see from the proofs that conditions of that
type suffice for some of our results.)
The condition we call \sdg\  is similar to the
conditions used in \cite{MP} and \cite{Gd}, and the
condition we call strict \sdg\  includes in addition a kind
of mixing condition on the summands.

This paper consists of four sections.
The first defines and proves useful relations between the
various forms of \sdg, and compares them with what is
already in the literature.
Section~2 contains the proofs of the first four parts
of Theorem~\ref{C01}.
These proofs follow by standard methods
from the work done in \cite{PhN}.
Section~3 contains the proof of the last part of
Theorem~\ref{C01}.
We were not able to follow the method of \cite{DNNP},
and in fact out proof does not use any version of the selection
theorem there.
Instead, we rely on perturbation results, functional
calculus, and a kind of approximate polar decomposition.
Finally, in Section~4 we give the applications to the
\ca s of \mh s, and the subalgebras of them considered in \cite{Ln}.

I am grateful to  Marius D\v{a}d\v{a}rlat, Anatole Katok,
Qing Lin, Cornel Pasnicu, and
Ian Putnam for useful conversations and email correspondence.
Most of this work was carried out during a sabbatical year at
Purdue University, and I am grateful to that institution
for its hospitality.

Some of the results of this paper were announced in \cite{LP}.

\section{Dimension growth}

The results on direct limits that we want to generalize from the
homogeneous case are mostly stated for systems with
slow dimension growth.
(See \cite{BDR}, \cite{Bl4}, \cite{Gd}, and \cite{MP}.)
We therefore discuss dimension growth in this section.
For direct systems of \rsha s, it is not clear what the appropriate
definition of slow dimension growth is.
(See the introduction for further discussion.)
We therefore confine ourselves to giving some usable definitions,
proving several easy results, and showing that our definitions are
satisfied for simple direct limits with no dimension growth.
These results suffice for our applications.

We state two versions of slow dimension growth.
The weaker version is more closely related to the definitions
in \cite{Gd} and \cite{MP}.
The stronger version includes a kind of mixing condition,
which seems to be needed in some of our proofs in the next section.

\begin{dfn}\label{G1}
Let $\left( \{A_i\}_{i \in \N}, \, \{ \ph_{i j} \} \right)$
be a direct system of
\rsha s, and let each $A_i$ be equipped with a specific decomposition
of length $l_i$ with total space $X_i$ and
topological dimension function $d_i \colon X_i \to \N \cup \{0\}$.
The system is said to have {\emph{slow dimension growth}}
(with respect to the given collection of decompositions)
if for every $i$, every \pj\  $p \in \Mi (A_i)$, and every $N \in \N$,
there is $j_0$ such that for all $j \geq j_0$ and $x \in X_j$ we have
\[
\ev_x (\ph_{i j} (p)) = 0 \,\,\,\,\,\,  {\text{or}} \,\,\,\,\,\,
   \rank ( \ev_x (\ph_{i j} (p))) \geq N d_j (x).
\]
The system is said to have {\emph{strict slow dimension growth}} if, in
the above, for $p \neq 0$ it is possible to choose $j_0$ such that
we always have $\rank ( \ev_x (\ph_{i j} (p))) \geq N d_j (x)$.
\end{dfn}

We note that $\rank ( \ev_x (\ph_{i j} (p)))$ depends only on the ranks
of the \pj s $\ev_y (p)$ for suitable $y \in X_i$, namely those $y$ for
which $\ev_y$ occurs among
the irreducible subrepresentations of the finite
dimensional representation $\ev_x \circ \ph_{i j}$ of $A_i$.
The global topological nature of $p$ is irrelevant.
The global topology does, however, have a strong influence on the
existence of \pj s $p$ with specified values of $\rank (\ev_y (p))$.
We would like to have $\rank (\ev_y (p)) = 1$.
However, it follows from Example~\ref{Q9} below that there are \rsha s
$A$ whose minimum matrix size is arbitrarily large but which contain
no nontrivial \pj s.
We do not even know whether, given $p \in \Mi (A)$, there is a
\pj\  $q \in A$ such that $\ev_y (q) = 0$ exactly when $\ev_y (p) = 0$.

This definition is complicated in practice, and we will therefore seek
simpler conditions which imply it.
First, we compare it with the definitions already in the literature
for the (nonsimple) homogeneous case.

\begin{prp}\label{G2}
Let $\left( \{A_i\}_{i \in \N}, \, \{ \ph_{i j} \} \right)$
be a direct system in which
\[
A_i = \bigoplus_{l = 1}^{r (i)}
            C \left( X_{i l}, \, M_{n (i, l)} \right),
\]
with $X_{i l}$ compact and connected.
Assume that the kernel of the map $A_i \to \dirlim A_j$
contains no entire summand $C \left( X_{i, k}, \, M_{n (i, k)} \right)$.
(That is, no summands vanish in the limit.)
Regard each $A_i$ as a \rsha\  with the obvious decomposition
(Example~1.4 of \cite{PhN}).
Then:
\bit
\item[(1)]
Slow dimension growth in the sense of \cite{Gd}, 2.1 implies slow
dimension growth in the sense of Definition~3.6 of \cite{MP}.
\item[(2)]
Slow dimension growth in the sense of \cite{Gd}, 2.1
implies slow dimension growth in the sense of Definition~\ref{G1}.
\item[(3)]
If the direct system has slow dimension growth in the
sense of Definition~3.6 of \cite{MP}, then it has a subsystem
which has slow dimension growth in the sense of Definition~\ref{G1}.
\item[(4)]
Slow dimension growth in the sense of Definition~\ref{G1}
implies slow
dimension growth in the sense of Definition~3.6 of \cite{MP}.
\eit
\end{prp}

We note that, for direct systems as in the proposition,
slow dimension growth in the sense of Definition~\ref{G1}
does not imply slow dimension growth in the sense of \cite{Gd}, 2.1.
Also, Definition~3.6 of \cite{MP} is formally stated only for
systems with unital injective maps, but makes sense in general.
The definition in \cite{Gd} is stated for direct systems over general
directed sets, but here we only consider direct systems over $\N$.

\medskip

{\emph{Proof of Proposition~\ref{G2}:}}
Let $e_i^{(l)}$ be the identity of the summand
$C \left( X_{i l}, \, M_{n (i, l)} \right)$.
Let
$\pi_i^{(l)} \colon A_i \to C \left( X_{i l}, \, M_{n (i, l)} \right)$
be the projection map.
For $j \geq i$ and $1 \leq m \leq r (j)$,
define the following quantities:
\[
\af_{i j}^{(m)}
  = \min \left\{ \rank {\ts{ \left( \rsz{ \pi_j^{(m)}
        \circ \ph_{i j} \left( \rsz{e_i^{(l)} } \right) } \right) }}
   \colon
   1 \leq l \leq r (i) \,\,{\text{and}} \,\,
   \pi_j^{(m)} \circ \ph_{i j} {\ts{ \left( \rsz{ e_i^{(l)} } \right) }}
       \neq 0 \right\},
\]
\[
\dt_j = \max \left\{ \frac{\dim (X_{j m})}{n (j, m)} \colon
      1 \leq m \leq r (j) \right\},
\]
and
\[
\mu_{i j} = \min  \left\{
   \frac{ \rank {\ts{ \left( \rsz{ \pi_j^{(m)} \circ \ph_{i j}
       \left( \rsz{e_i^{(l)} } \right) } \right) }}}{n (j, m)} \colon
   1 \leq l \leq r (i), \,\, 1 \leq m \leq r (j), \,\,{\text{and}} \,\,
   \pi_j^{(m)} \circ \ph_{i j} \left( \rsz{ e_i^{(l)} } \right) \neq 0
            \right\}.
\]
Note that the ranks appearing here are constant functions,
which we identify with the corresponding integers, because the
spaces $X_{j m}$ are connected.
The quantity $\dt_j$ is called $d_j$ in \cite{Gd}, 2.1.
Also, $\ph_{i j}$ and $\mu_{i j}$ are written as
$\ph_{j i}$ and $\mu_{j i}$ there.
With these definitions, our system has
slow dimension growth in the sense of \cite{Gd}, 2.1
\ifo\  for every $i \in \N$ we have
\[
\lim_{j \to \infty} \frac{\dt_j}{\mu_{i j}} = 0,
\]
and it has
slow dimension growth in the sense of Definition~3.6 of \cite{MP}
\ifo\  for every $i \in \N$ we have
\[
\liminf_{j \to \infty}
   \max  \left\{ \frac{\dim (X_{j m})}{\af_{i j}^{(m)} } \colon
         1 \leq m \leq r (j) \right\} = 0.
\]
(Here, and in the rest of this proof,
such limits are taken over $j \geq i$.)

We now claim that a system has
slow dimension growth in the sense of Definition~\ref{G1}
\ifo\  for every $i \in N$ we have
\[
\lim_{j \to \infty}
     \max  \left\{ \frac{\dim (X_{j m})}{\af_{i j}^{(m)} } \colon
         1 \leq m \leq r (j) \right\} = 0.
\]

To see this, first assume the condition of Definition~\ref{G1} holds.
Choose $i$ and $N$, and apply the condition to the \pj s $e_i^{(l)}$ for
$1 \leq l \leq r (i)$.
Call the resulting numbers $j_0 (l)$.
Set $j_0 = \max \{ j_0 (l) \colon 1 \leq l \leq r (i) \}$.
Then for all $j \geq j_0$ and all $m$ with
$1 \leq m \leq r (j)$, we have
\[
\frac{\dim (X_{j m})}{\af_{i j}^{(m)} } \leq \frac{1}{N}.
\]
Thus, the limit condition above is satisfied.

Conversely, assume the limit condition above, and let $p \in \Mi (A_i)$
be a \pj.
For $1 \leq l \leq r (i)$ set
$\rh_l = \rank {\ts{ \left( \rsz{ \pi_i^{(l)} (p) } \right) }} /
      \rank {\ts{ \left( \rsz{ e_i^{(l)} } \right) }}$.
Set
\[
\rh =
\min \{ \rh_l \colon 1 \leq l \leq r (i)
                 \,\, {\text{and}} \,\, \rh_l \neq 0\}.
\]
Note that the \pj\   $\pi_j^{(m)} \circ \ph_{i j} (p)$
has constant rank equal to
\[
\sum_{l = 1}^{r (i)}
    \rh_m \rank {\ts{ \left( \rsz{ \pi_j^{(m)} \circ \ph_{i j}
         {\ts{ \left( \rsz{ e_i^{(l)} } \right) }} } \right) }}.
\]
If $\pi_j^{(m)} \circ \ph_{i j} (p) \neq 0$, then there is $l_0$ with
$\rh_{l_0} \neq 0$ and
$\rank \left( \rsz{ \pi_j^{(m)} \circ \ph_{i j}
        \left( \rsz{ e_i^{(l_0)} } \right) } \right) \neq 0$.
So
\[
\rank {\ts{ \left( \rsz{ \pi_j^{(m)} \circ \ph_{i j} (p)  } \right) }}
    \geq \rh_{l_0}
        \rank {\ts{ \left( \rsz{ \pi_j^{(m)} \circ \ph_{i j}
           \left( \rsz{ e_i^{(l_0)} } \right) } \right) }}
    \geq \rh \af_{i j}^{(m)}.
\]
According to the limit condition above, we can choose $j_0$ such that
if $j \geq j_0$ and $1 \leq m \leq r (j)$, then
\[
\frac{ \dim (X_{j m}) }{ \af_{i j}^{(m)} } \leq \frac{ \rh}{N}.
\]
For such $j$, and whenever $\pi_j^{(m)} \circ \ph_{i j} (p) \neq 0$, we
have
\[
\rank {\ts{ \left( \rsz{ \pi_j^{(m)} \circ \ph_{i j} (p)  } \right) }}
        \geq N \dim (X_{j m}).
\]
This verifies slow dimension growth in the sense of Definition~\ref{G1},
and proves the claim.

Part~(4) of the proposition is immediate from the claim.
To get part~(2), we merely observe that for $1 \leq m \leq r (j)$
we have
\[
\frac{ \dim (X_{j m}) }{\af_{i j}^{(m)} } \leq \frac{\dt_j}{\mu_{i j} }.
\]
Part~(1) follows immediately from parts~(2) and~(4).

It remains to prove part~(3).
Given the claim, we must prove that we can replace the $\liminf$ in
Definition~3.6 of \cite{MP} by a limit by passing to a cofinal subset.
We do this by a kind of diagonalization argument.
So assume that for all $i$ we have
\[
\liminf_{j \to \infty}
    \max  \left\{ \frac{\dim (X_{j m})}{\af_{i j}^{(m)} } \colon
         1 \leq m \leq r (j) \right\} = 0.
\]
Set $i (1) = 1$, and choose an infinite subset $I_1 \subset \N$,
with $i (1) \in I_1$, such that
\[
\lim_{j \to \infty, \, j \in I_1}
    \max  \left\{ \frac{\dim (X_{j m})}{\af_{i (1) j}^{(m)} } \colon
         1 \leq m \leq r (j) \right\} = 0.
\]
Let $i (2)$ be the second element of $I_1$, and choose an infinite
subset $I_2 \subset I_1$, with $i (1), i (2) \in I_2$, such that
\[
\lim_{j \to \infty, \, j \in I_2}
    \max  \left\{ \frac{\dim (X_{j m})}{\af_{i (2) j}^{(m)} } \colon
         1 \leq m \leq r (j) \right\} = 0.
\]
Proceed inductively.
Take
\[
I = \bigcap_{n = 1}^{\infty} I_n = \{ i (1), i(2), \dots \},
\]
which is cofinal in $\N$.
Then
\[
\lim_{j \to \infty, \, j \in I}
    \max  \left\{ \frac{\dim (X_{j m})}{\af_{i j}^{(m)} } \colon
         1 \leq m \leq r (j) \right\} = 0
\]
for all $i \in I$, so that the corresponding subsequence has
slow dimension growth in the sense of Definition~\ref{G1}.
\QED

\medskip

We can generalize the approach of \cite{Gd} and \cite{MP} somewhat.

\begin{dfn}\label{G3}
Let $\ph \colon A \to B$ be a (not necessarily unital) \hm\  of \rsha s.
Let decompositions of $A$ over $X = \coprod_{k = 0}^K X_k$ and
$B$ over $Y = \coprod_{l = 0}^L Y_l$ be given, with all $X_k$ and
$Y_l$ connected.
For $y \in Y$ define $\mu_{k, y} (\ph)$, the
{\emph{$k$-th partial multiplicity of $\ph$ at $y$,}} as follows.
Consider all possible direct sum decompositions
$\ev_y \circ \ph \cong 0 \oplus \bigoplus_{j = 1}^n \ev_{x_j}$
with $x_j \in X$.
(Because $\ev_y$ is a finite dimensional representation, Lemma~2.1
of \cite{PhN}
implies that there is always at least one such decomposition.)
Then $\mu_{k, y} (\ph)$ is the maximum, over all such
decompositions, of the number (counting multiplicity)
of $x_j$ that are in $X_k$.
Moreover, define $\mu_{k, l} (\ph)$, the
{\emph{$k$-th partial multiplicity of $\ph$ at $Y_l$,}} to be
$\sup_{y \in Y_l} \mu_{k, y} (\ph)$.
Finally, say that $\ph$ is
{\emph{zero in the $(k, l)$-component}} if for all $a \in A$
there is $\widetilde{a} \in A$ such that
$\ev_x (\widetilde{a}) = \ev_x (a)$ for $x \in  X_k$ and
$\ev_y ( \ph (\widetilde{a})) = 0$ for all $y \in Y_l$.
\end{dfn}

Note that, for direct sums of trivial homogeneous \ca s,
zero in the $(k, l)$-component
simply means that the partial map from the
$k$-th summand of $A$ to the $l$-th summand of $B$ is zero.

\begin{lem}\label{G4}
Let $\left( \{A_i\}_{i \in \N}, \, \{ \ph_{i j} \} \right)$
be a direct system of
\rsha s, and let each $A_i$ be equipped with a specific decomposition
of length $L_i$, with {\emph{connected}} base spaces
$X_{1, 0}, X_{i, 1}, \dots, X_{i, L_i}$
and with total space $X_i = \coprod_{k = 0}^{L_i} X_{i, k}$.
Assume that for every $i$ and every $N \in \N$,
there is $j_0$ such that,
for all $j \geq j_0$, all $k$ with $0 \leq k \leq L_i$, and all
$l$ with $0 \leq l \leq L_j$, either
$\ph_{i j}$ is zero in the $(k, l)$-component or
$\mu_{k, l} (\ph_{i j}) \geq N \dim (X_{j, k})$.
Then the direct system has slow dimension growth.
\end{lem}

\begin{pff}
Let $p \in \Mi (A_i)$,
choose $j_0$ as in the hypotheses of the lemma, and
let $j \geq j_0$.
Write $p = (p_0, p_1, \dots, p_{L_i})$, where $p_k$ is the restriction
of the standard representation of $p$ (a function on $X_i$) to
$X_{i, k}$.
Let $q = \ph_{i j} (p)$, and analogously write
$q = (q_0, q_1, \dots, q_{L_j})$.
Since the spaces are all connected, the $p_k$ and $q_l$ have constant
ranks.

Let $0 \leq l \leq L_j$.
Suppose that there is some $k$ with $p_k \neq 0$ and
$\mu_{k, l} (\ph_{i j}) \geq N \dim (X_{j, l})$.
Choose $y \in X_{k, l}$ such that
$\mu_{k, y} (\ph_{i j}) \geq N \dim (X_{j, l})$.
Then clearly
\[
\rank (\ev_y (q)) \geq
             N \dim (X_{j, l}) \rank (p_k) \geq N \dim (X_{j, l}).
\]
Since $q_l$ has constant rank, this inequality holds for arbitrary
$y \in X_{k, l}$.

Otherwise, $\ph_{i j}$ is zero in the $(k, l)$-component for all $k$
with $p_k \neq 0$.
For such $k$, choose $a_k \in \Mi (A_i)$ such that
$\ev_x (a_k) = \ev_x (p)$ for $x \in  X_{k, i}$ and
$\ev_y ( \ph_{i j} (a_k)) = 0$ for all $y \in X_{j, l}$.
(This can clearly be done by considering the entries separately.)
Set $a = \sum_{k\colon p_k \neq 0} a_k^* a_k$.
Then $a \geq p$ and $\ev_y ( \ph_{i j} (a)) = 0$
for all $y \in X_{j, l}$.
Therefore $q_l = 0$.
\end{pff}

\medskip

For present applications, we are primarily interested in the simple
case, and we devote the rest of this section to it.
In this case, at least with no (rather than slow) dimension growth,
we do not need connectedness assumptions.

The following lemma is a slight generalization of part of
Proposition~2.1 of \cite{DNNP}.

\begin{lem}\label{G5}
Let $A = \dirlim A_i$ be a simple direct limit of \rsha s, such that
all the maps $\ph_{i j} \colon A_i \to A_j$ in the system are unital and
injective.
Let $X_i$ be the total space of $A_i$.
Let $a \in A_i \setminus \{ 0 \}$ for some $i$.
Then there exists $j_0$ such that, for every $j \geq j_0$ and every
$x \in X_j$, we have $\ev_x ( \ph_{i j} (a)) \neq 0$.
\end{lem}

\begin{pff}
\Wolog\ $i = 0$.
Assume the conclusion fails for some $a$.
Passing to a subsequence in the direct system, we may assume
that for every $j$ there is $x \in X_j$ such that
$\ev_x ( \ph_{0j} (a)) = 0$.
Therefore the ideal $I_j = \overline{A_j \ph_{0j} (a) A_j}$
is nontrivial.
Since $\ph_{j, j + 1} (I_j) \subset I_{j + 1}$, we may form the
ideal $I = \dirlim I_j \subset A$.
This ideal is nonzero since it contains the image of $a$.
If $1 \in I$, then (using injectivity of the $\ph_{i j}$) there
is $j$ and $b \in I_j$ with $\| b - 1 \| < \frac{1}{2}$,
which contradicts $I_j \neq A_j$.
So $I$ is a proper ideal in $A$, contradicting simplicity.
\end{pff}

\begin{lem}\label{G7}
Let $A$ be any \ca\  in which there do not exist $n + 1$ mutually
orthogonal nonzero selfadjoint elements.
Then $\dim (A) \leq n^2$.
\end{lem}

\begin{pff}
\Wolog\  there are $n$ mutually orthogonal nonzero selfadjoint elements
$a_1, \dots, a_n \in A$.
Then there is a unique nonzero $\af_j \in \spec (a_j)$ for each $j$.
(If some $\spec (a_j)$ has more than one nonzero element, then
\ct\  functional calculus gives two
mutually orthogonal nonzero selfadjoint elements $b, \, c \leq | a_j|$,
and using them in place of $a_j$ contradicts the assumption.)
The elements $p_j = \af_j^{-1} a_j$ are mutually orthogonal nonzero
\pj s.
Moreover, $p = \sum_{j = 1}^n p_j$ must be an identity for $A$,
since the existence of a nonzero selfadjoint element of
$(1 - p) A (1 - p)$ contradicts the assumption.

If $p_j A p_j$ contains a selfadjoint element $b$ not a scalar multiple
of $p_j$, then $\spec (b)$
(taken relative to the unital algebra $p_j A p_j$)
has at least two elements.
We can then use \ct\  functional calculus to get a contradiction
as in the previous paragraph.
So $p_j A p_j = \C \cdot p_j$ for all $j$.

Now suppose $p_j A p_k \neq 0$ for some $j$ and $k$.
Let $c \in p_j A p_k$ be nonzero.
Then there are $\ld, \, \mu \in (0, \infty)$ such that
$c c^* = \ld p_j$ and $c^* c = \mu p_k$.
Replacing $c$ by a suitable scalar multiple, we can assume $\ld = 1$.
Then also $\mu = 1$.
Let $d \in p_j A p_k$ be arbitrary.
Then $c^* d \in p_k A p_k$, so $c^* d = \af p_k$ for some $\af \in \C$.
It follows that
\[
d = p_j d = c c^* d = c \cdot  \af p_k = \af c.
\]
This computation shows that $\dim (p_j A p_k) \leq 1$.
Since $A = \bigoplus_{j,\, k} p_j A p_k$ as a Banach space,
we conclude that $\dim (A) \leq n^2$.
\end{pff}

\begin{lem}\label{G7.5}
Let $A$ be a simple \ca.
Suppose $A$ has a hereditary subalgebra $B$ such that $B \cong M_n$.
Then $A$ is isomorphic to the algebra $K (H)$
of all compact operators on some Hilbert space $H$.
\end{lem}

\begin{pff}
\Wolog\ $B \cong \C$, that is, $B = \C \cdot p$ for some
\pj\  $p \in A$.
Then there is a state $\om$ on $A$ such that $p a p = \om (a) p$
for all $a \in A$.
With the help of this state, it is easy to make $p A$ a
$\C$--$A$ strong Morita equivalence bimodule.
(It is full as an $A$-module because $A$ is simple.)
In particular, $H = p A$ is a Hilbert space such that $A \cong K (H)$.
\end{pff}

\begin{lem}\label{G8}
Let $A = \dirlim A_i$ be an infinite dimensional
simple direct limit of \rsha s, such that
all the maps $\ph_{i j} \colon A_i \to A_j$ in the system are unital and
injective.
Let $X_i$ be the total space of $A_i$.
Let $a \in A_i \setminus \{ 0 \}$ for some $i$. 
Then for every $n \in \N$ there exists $j_0$ such that,
for every $j \geq j_0$ and every
$x \in X_j$, we have
$\rank \left( \ev_x ( \ph_{i j} (a)) \right) \geq n$.
\end{lem}

\begin{pff}
\Wolog\  $i = 0$.
Moreover, since $\rank (b^* b) = \rank (b)$ for any $b \in M_n$,
\wolog\  $a \geq 0$.

Choose $l$ such that there are $n$
mutually orthogonal nonzero selfadjoint elements
$b_1, \dots, b_n \in \overline{\ph_{0 l} (a) A_{l} \ph_{0 l} (a)}$.
(If this is not possible, then Lemma~\ref{G7} implies that
$\dim \left( \rsz{ \overline{\ph_{0 j} (a) A_j \ph_{0 j} (a)} } \right)
                             \leq (n - 1)^2$
for all $j$.
Then, with $c$ being the image of $a$ in $A$, we have
$\dim ( \overline{c A c}) \leq (n - 1)^2$.
Since $A$ is simple and unital, Lemma~\ref{G7.5} implies that
$A \cong M_m$ for some $m$.
This contradicts infinite dimensionality.)
Now choose (by Lemma~\ref{G5}) $j_0 \geq l$  such that,
for every $j \geq j_0$, every $k$, and every $x \in X_j$, we have
$\ev_x ( \ph_{l j} (b_k)) \neq 0$.
For such $j$, and $x \in X_j$, the image
\[
\ev_x
    \left( \rsz{ \overline{\ph_{0 j} (a) A_{j} \ph_{0 j} (a)} } \right)
  =  \overline{ \ev_x ( \ph_{0 j} (a)) M_m \ev_x ( \ph_{0 j} (a))}
\]
(for some suitable $m$) contains $n$
mutually orthogonal nonzero selfadjoint elements
$\ev_x (b_1), \dots, \ev_x (b_n)$.
Therefore $\rank \left( \ev_x ( \ph_{0 j} (a)) \right) \geq n$.
\end{pff}

\begin{cor}\label{G9}
Let $A = \dirlim A_i$ be an infinite dimensional
simple direct limit of \rsha s, such that
all the maps $\ph_{i j} \colon A_i \to A_j$ of the system are unital and
injective.
Assume the system has no dimension growth, that is, there is
a finite $d$ such that all $A_i$ have \tdim\  at most $d$.
Then the system has strict slow dimension growth.
\end{cor}

\begin{sloppypar}
\begin{pff}
We verify Definition~\ref{G1} for a particular \pj\  $p$ by
applying Lemma~\ref{G8} to the algebras $M_n (A) = \dirlim M_n (A_i)$,
taking $a = p$.
\end{pff}
\end{sloppypar}

\medskip

We now briefly consider the effect on dimension growth of forcing the
maps of a system to be injective.
This means that the algebras in the system must be replaced by
quotients.
Proposition~3.1 of \cite{PhN}
shows that this does not increase the \tdim.
However, (strict) slow dimension growth also depends on the details
of the decompositions, and we do not know how to obtain a decomposition
of a quotient which is suitably related to a given
decomposition of the original algebra.
For simple direct limits with no dimension growth, the special case
we are most interested in, the previous corollary eliminates this
difficulty.

\begin{prp}\label{G10}
Let $A = \dirlim A_i$ be a simple direct limit of separable
\rsha s which has
no dimension growth in the sense of Corollary~\ref{G9}.
Then $A$ is the direct limit of a direct system $\{B_i\}$ with
no dimension growth, and in which in addition all the maps
$\ps_{i} \colon B_i \to B_{i + 1}$ of the system are injective.
\end{prp}

\begin{pff}
Let $d$ be a finite upper bound on the \tdim s  of the $A_i$.
Let $B_i$ be the image of $A_i$ in $A$.
Then also $\dirlim B_i \cong A$.
Since the $A_i$ are separable,
Proposition~3.1 of \cite{PhN}
implies that the $B_i$ are separable \rsha s
with \tdim\  at most $d$.
\end{pff}

\section{Cancellation and comparison in direct limits}

In this section, we prove the positive results on direct limits with
slow dimension growth, except for stable rank $1$.
The first result (essentially cancellation for unitaries) seems
not to have been noticed before, but the others are all analogs of
known results in the homogeneous case.
We don't know whether Theorem~\ref{C2} (cancellation) remains true
without strict slow dimension growth, or whether Theorem~\ref{C4} (weak
unperforation) remains true without simplicity.

We impose injectivity on the maps of the system when convenient.
Proposition~\ref{G10} shows this condition can be eliminated for
simple direct limits with no dimension growth.

\begin{thm}\label{C1}
Let $A = \dirlim (A_i, \ph_{i j})$ be a unital direct limit
of a system of \rsha s with slow dimension growth.
Then the map $U (A) / U_0 (A) \to K_1 (A)$ is an isomorphism.
\end{thm}

\begin{pff}
\Wolog\  all maps of the system are unital.
Let $\ph_{i \infty} \colon A_i \to A$ be the induced maps to
the direct limit; these are also unital.

We first prove surjectivity.
Let $\et \in K_1 (A)$.
Choose $i$, $n$, and $u \in U (M_n (A_i))$ such that
$(\ph_{i \infty})_* ([u]) = \et$.
Let $q$ be the identity of $M_n (A_i)$, and let $p$
be the identity of $A_i$, regarded as a subalgebra of $M_n (A_i)$
by identifying it with the upper left corner.
Apply the definition of slow dimension growth to the \pj\  $p$.
This gives $j \geq i$ such that, if $X$ is the total space of $A_j$
and $d$ is its \tdim\  function, then
$\rank (\ev_x (\ph_{i j} (p))) \geq \frac{1}{2} d (x)$ for all $x \in X$.
(Note that $\ev_x (\ph_{i j} (p))$ is never zero, because
$\ph_{i j} (p)$ is the identity of $A_j$.)

Apply Proposition~4.4~(1) of \cite{PhN}
to $\ph_{i j} (p)$, $\ph_{i j} (q)$, and
$\ph_{i j} (u)$,
obtaining $v \in U (A_j)$ such that $v \oplus 1$ is homotopic to
$\ph_{i j} (u)$ in $U (M_n (A_j))$.
Then $\ph_{j \infty} (v) \in U (A)$ and $[\ph_{j \infty} (v)] = \et$
in $K_1 (A)$.

Now we prove injectivity.
Let $u \in U (A)$ satisfy $[u] = 0$ in $K_1 (A)$.
We show $u \in U_0 (A)$.
We have $u \oplus 1 \in U_0 (M_n (A))$ for some $n$.
By standard methods there is $i$ and $u_0 \in U (A_i)$
such that $\ph_{i \infty} (u_0)$ is homotopic to $u$ in $U (A)$
and such that $u_0 \oplus 1 \in U_0 (M_n (A_i))$.
It suffices to find $j \geq i$ such that $\ph_{i j} (u) \in U_0 (A_j)$.

As in the proof of surjectivity, let $q$ be the identity of
$M_n (A_i)$, and let $p$ be the identity of $A_i$.
For the same reason as there, we can choose
$j \geq i$ such that, with $X$ and $d$ as there, we have
$\rank (\ev_x (\ph_{i j} (p))) \geq d (x)$ for all $x \in X$.
Since also $\rank (\ev_x (\ph_{i j} (p))) \geq 1$ for all $x \in X$,
it follows that
$\rank (\ev_x (\ph_{i j} (p))) \geq \frac{1}{2} (d (x) + 1)$
for all $x \in X$.
Using $\ph_{i j} (u_0) \oplus 1 \in U_0 (M_n (A_j))$,
apply Proposition~4.4~(2) of \cite{PhN}
to obtain $\ph_{i j} (u_0) \in U_0 (A_j)$, as desired.
\end{pff}

\medskip

The next theorem is the analog in our situation of Proposition~3.7~(a)
of \cite{Bl4}, of one part of Theorem~3.7 of \cite{MP},
and of Proposition~2.6 of \cite{Gd}.

\begin{thm}\label{C2}
Let $A = \dirlim (A_i, \ph_{i j})$ be a direct limit
of a system of \rsha s with strict slow dimension growth.
Then the \pj s in $\Mi (A)$ satisfy cancellation:
if $e \oplus q \sim f \oplus q$, then $e \sim f$.
\end{thm}

\begin{pff}
By standard arguments, we may assume that $e$, $f$, and $q$ are
in $\Mi (A_i)$ for some $i$, and that
$e \oplus q \sim f \oplus q$ in $\Mi (A_i)$.
Suppose first that $e = 0$.
Let $X$ be the total space of $A_i$.
For every $x \in X$ we have $\ev_x (q) \sim \ev_x (f) \oplus \ev_x (q)$,
whence $f = 0$.
So certainly $e \sim f$.
Otherwise, use strict slow dimension growth to choose $j \geq i$
such that, if $X$ is the total space of $A_j$
and $d$ is its \tdim\  function, then
$\rank (\ev_x (\ph_{i j} (e))) \geq \frac{1}{2} d (x)$ for all $x \in X$.
Under these conditions,
$\ph_{i j} (e) \sim \ph_{i j} (f)$ by Proposition~4.3~(2) of \cite{PhN}.
So $e \sim f$ in $\Mi (A_i)$.
\end{pff}

\medskip

The next theorem is the analog in our situation of
the other part of Theorem~3.7 of \cite{MP}
and of one part of Proposition~3.7~(b) of \cite{Bl4}.

\begin{thm}\label{C3}
Let $A = \dirlim (A_i, \ph_{i j})$ be a direct limit of a system of
\rsha s with slow dimension growth, and in which the maps $\ph_{i j}$
of the system are all injective and unital.
Then $A$ satisfies
Blackadar's Second Fundamental Comparability Question
(\cite{Bl3}, 1.3.1):
if $p, \, q \in \Mi (A)$ are \pj s such that $\ta (p) < \ta (q)$
for all normalized traces $\ta$ on $A$, then $p \precsim q$.
\end{thm}

\begin{pff}
By standard arguments we may assume that $p$ and $q$ are in $M_n (A_i)$
for some $n$ and $i$.
Dropping initial terms, \wolog\  $i = 0$.
Replacing every $A_i$ by $M_n (A_i)$ does not change any of the
hypotheses (see Lemma~1.12 of \cite{PhN}),
so \wolog\  $p, \, q \in A_0$.
Arguing as in the proof of Proposition~4.1 of \cite{Bl1},
we find $i$ such that
$\ta (\ph_{0 i} (p)) < \ta (\ph_{0 i} (q))$
for all normalized traces $\ta$ on $A_i$.
Letting $X$ be the total space of $A_i$, letting
$m \colon X \to \N \cup \{0\}$ be its matrix size function
(as in Definition~\ref{A2}),
letting $\Tr$ be the usual trace on matrices, and
taking $\ta = \frac{1}{m (x)} \Tr \circ \ev_x$, we obtain
$\rank (\ev_x (\ph_{0 i} (p))) < \rank (\ev_x (\ph_{0 i} (q)))$
for all $x \in X$.
Let $N = \max_{x \in X} \rank (\ev_x (\ph_{0 i} (q)))$.
Note that $N > 0$.
The slow dimension growth hypothesis provides
$j \geq i$ such that, with
$Y$ being the total space of $A_j$ and
$d \colon Y \to \N \cup \{0\}$ being its topological dimension function
(as in Definition~\ref{A2}),
we have
\[
\ev_y (\ph_{0 j} (q)) = 0 \,\,\,\,\,\,  {\text{or}} \,\,\,\,\,\,
   \rank ( \ev_y (\ph_{0 j} (q))) \geq N d (y)
\]
for all $y \in Y$.
For $y$ for which the second case holds,
note that $\ev_y \circ \ph_{i j}$ is a finite dimensional
representation of $A_i$, and therefore is equivalent to a direct sum
$0 \oplus \bigoplus_{r = 1}^R \ev_{x_r}$ with $x_1, \dots, x_R \in X$.
Since
\[
\rank ( \ev_y (\ph_{0 j} (q))) =
                       \sum_{r = 1}^R \rank ( \ev_{x_r} (\ph_{0 i} (q)))
            \andeqn
\rank (\ev_{x_r}  (\ph_{0 i} (q))) \leq N,
\]
we get $R \geq d (y)$.
We have
$\rank (\ev_x (\ph_{0 i} (q))) \geq \rank (\ev_x (\ph_{0 i} (p))) + 1$
for all $x \in X$, so
\[
\rank ( \ev_y (\ph_{0 j} (q)))
              \geq \rank ( \ev_y (\ph_{0 j} (p))) + d (y).
\]
It follows from Proposition~4.3~(1) of \cite{PhN}
that $\ph_{0 j} (p) \precsim \ph_{0j} (q)$ in $A_j$.
This proves the result.
\end{pff}

\medskip

Finally, we deal with unperforation.
At least two different definitions of weak unperforation appear in the
literature, namely Definition~6.7.1 of \cite{Bl2} and 2.1 of \cite{El1}.
These agree in the simple case (2.1 of \cite{El1}), which is the one
relevant below.
To minimize confusion, however, we use the more descriptive term
``unperforated for the strict order'' for the property we prove
(Definition~6.7.1 of \cite{Bl2}).

The following result is the analog of another part of
Proposition~3.7~(b) of \cite{Bl4}, and of Proposition~2.10 of \cite{Gd}.
Unfortunately, we have not been able to generalize the method
of \cite{Gd}, so we get the result only in the simple case.
As observed in \cite{Bl4}, in that case it is immediate from
the Second Fundamental Comparability Question.

\begin{thm}\label{C4}
Let $A = \dirlim (A_i, \ph_{i j})$ be a simple direct limit of a system
of \rsha s with slow dimension growth, and in which the maps $\ph_{i j}$
of the system are all injective and unital.
Then $K_0 (A)$ is unperforated for the strict order.
That is, if $\et \in K_0 (A)$ and there is $n > 0$ such that
$n \et > 0$, then $\et > 0$.
\end{thm}

\begin{pff}
Write $\et = [p] - [q]$ for \pj s $p, \, q \in \Mi (A)$.
Since $A$ is simple and
$p \oplus \cdots \oplus p \sim q \oplus \cdots \oplus q \oplus r$
for some nonzero \pj\  $r$ ($n$ copies each of $p$ and $q$),
we have $\ta (p) > \ta (q)$ for all normalized traces $\ta$ on $A$.
So $p \succsim q$ by Theorem~\ref{C3}.
\end{pff}

\section{Reduction of stable rank}

In this section, we prove that if $A$ is a simple direct limit of a
system of
\rsha s with no dimension growth, then $A$ has stable rank $1$.
This generalizes the result of \cite{DNNP}.
Our proof, however, is somewhat different, being based on a
notion of an ``approximate polar decomposition''.
Essentially, if $a$ is an element of a \rsha,
and if each $\ev_x (a)$ is small on a sufficiently large subspace
of the (finite dimensional) space on which it acts,
then there is a unitary $u$
such that $u (a^* a)^{1/2}$ is close to $a$.

We start with two preparatory results.
The first, which provides a means of constructing continuous
projection valued functions without appealing to any
selection theorems,
has independent usefulness.
Most of the work is contained in the third lemma,
which essentially does the induction step in the construction
of the approximate polar decomposition by induction on the length
of a \rshd.

\begin{prp}\label{C9}
Let $X$ be a \cpt\  space, let $a \in C (X, M_n)_{\sa}$, and let
$\ld \in \R$.
For $x \in X$ define a projection $p (x)$ by
$p (x) = \ch_{( - \infty, \ld)} (a (x))$.
Then there exist open sets $U_k$, for $0 \leq k \leq n$,
and \ct\  rank $k$ \pj s $p_k \colon U_k \to M_n$, such that:
\bit
\item[(1)]
$\bigcup_{k = 0}^n U_k = X$.
\item[(2)]
If $k \leq l$ and $x \in U_k \cap U_l$, then $p_k (x) \leq p_l (x)$.
\item[(3)]
$p (x) = \sup \{ p_k (x) \colon x \in U_k\}$ for all $x \in X$.
\item[(4)]
$p_k (x)$ commutes with $a (x)$ for all $x \in U_k$.
\eit
\end{prp}

\begin{pff}
For $x \in X$, write the eigenvalues of $a (x)$ as
\[
\af_1 (x) \leq \af_2 (x) \leq \cdots \leq \af_n (x)
\]
(repeated according to multiplicity).
It follows from Theorem~8.1 of \cite{Bh} that the $\af_k$ are continuous
functions on $X$.
Further set $\af_{n + 1} (x) = \infty$ for all $x$.
Define
\[
\bt_k (x) = {\textstyle{\frac{1}{2}}} \left( \af_k (x) + \ld \right)
            \andeqn
U_k = \{ x \in X \colon \af_k (x) < \bt_k (x) < \af_{k + 1} (x) \}
\]
for $1 \leq k \leq n$.
Then define $p_k (x) = \ch_{( - \infty, \bt_k (x))} (a (x))$ for
$x \in U_k$.
Further take $U_0 = X$ and $p_0 (x) = 0$.

We verify that these sets and projections satisfy the conclusion of the
proposition.
The $U_k$ are open because the functions $\af_k$ and $\bt_k$ are \ct.
To see that $p_k$ is \ct, rewrite $p_k (x) = f_x (a (x))$, where
\[
f_x (t) = \left\{ \begin{array}{lcl}
  1 &  &  t \leq \af_k (x) \\
       {\displaystyle{
               \frac{\af_{k + 1} (x) - t}{\af_{k + 1} (x) - \af_k (x)}}}
   &   &  \af_k (x) \leq t \leq \af_{k + 1} (x)  \\
  0 &  &  \af_{k + 1} (x) \leq t.
         \end{array} \right.
\]
The function $(t, x) \mapsto f_x (t)$ is jointly \ct,
so $x \mapsto f_x (a (x))$ is \ct\  by Proposition~2.12 of \cite{Ph0}.
Clearly $\rank (p_k (x)) = k$ for all $x$.
It is obvious that the $U_k$ cover $X$.
It is also obvious that if $k \leq l$ then $p_k (x) \leq p_l (x)$
wherever both are defined.
To verify part~(3), we note that if $x \in U_k$, then
${\textstyle{\frac{1}{2}}} \left( \af_k (x) + \ld \right) > \af_k (x)$,
whence $\af_k (x) < \ld$, so that $p_k (x) \leq p (x)$.
On the other hand, for $x \in X$ there is some $k$ with
$\af_k (x) < \ld \leq \af_{k + 1} (x)$.
For this $k$, we have $x \in U_k$ and $p_k (x) = p (x)$.
We have thus proved that
$p (x) = \sup \{ p_k (x) \colon x \in U_k\}$.
Finally, part~(4) is immediate because
$p_k (x)$ is obtained from $a (x)$ using functional calculus.
\end{pff}

\begin{lem}\label{D1}
Let $\ep > 0$.
Then there is $\dt > 0$ such that the following holds.
Let $Z = Z_1 \cup Z_2$ be a \cpt\  space, with $Z_1$ and $Z_2$
closed subsets, and let $U$ be a \nbhd\  of $Z_1$.
Let $A$ be a unital \ca, and let $p \in C(Z, A)$ be a \pj.
Let $s_j \in C (Z_j, A)$ be a partial isometry with
initial \pj\  $s_j^* s_j = p |_{Z_j}$ for $j = 1, \, 2$,
and suppose $\| s_1 |_{Z_1 \cap Z_2} - s_2 |_{Z_1 \cap Z_2} \| < \dt$.
Then there is a partial isometry $s \in C (Z, A)$ with
initial \pj\  $s^* s= p$, such that
\[
s |_{Z_1} = s_1, \,\,\,\,\,\,
   s |_{Z \setminus U} = s_2  |_{Z \setminus U},
  \andeqn \left\| s |_{Z_2} - s_2 \right\| < \ep.
\]
\end{lem}

\begin{pff}
The partial isometry $s$ will be constructed as follows.
Choose some \ct\  function $c \colon Z \to A$ such that $c |_{Z_1} = s_1$.
Let
\[
V = Z_1 \cup \{x \in Z_2 \colon \| c (x) p (x) - s_2 (x) \| < \dt \},
\]
which is a \nbhd\  of $Z_1$.
Choose a \ct\  function $f \colon Z \to [0, 1]$ which is equal to
$1$ on $Z_1$ and equal to zero on $Z \setminus (U \cap V)$.
Then define
\[
a (x) = f (x) c (x) p (x) + (1 - f (x)) s_2 (x)  \andeqn
  s (x) = a (x)  \left[a (x)^* a (x) \right]^{-1/2},
\]
with functional calculus in $p (x) A p (x)$, noting that
$a (x) p (x) = p (x)$ for all $x$.
It is clear that if $\| c (x) p (x) - s_2 (x)\|$ is small enough
for $x \notin U \cap V$
(depending only on $\ep$), then we will get
$\| s (x) - s_2 (x)\| < \ep$ for $x \in Z_2$.
\end{pff}

\medskip

The following lemma is the heart of the construction of
the approximate polar decomposition.
It is a relative version of the result for $C (X, M_n)$.

\begin{lem}\label{D2}
Let $\af, \, \ep > 0$ and $n \in \N$.
Then there is $\dt > 0$ such that the following holds.
Let $X$ be a \cpt\  space with dimension $d = \dim (X)$, and
let $X^{(0)} \subset X$ be closed.
Let $a \in C (X, M_n)$ satisfy $\| a \| \leq 1$,
and assume that for each $x \in X$
there is a subspace $E_x$ of $\C^n$ with $\dim (E_x) \geq \frac{1}{2} d$
such that
$\| a (x) \xi \| < \af \| \xi \|$ for $\xi \in E_x \setminus \{ 0\}$.
Let $p$ be the \lscp\  %
\[
x \mapsto \ch_{ ( -\infty, \af)}
   {\ts{ \left( \rsz{ \left[a (x)^* a (x) \right]^{1/2} } \right) }}.
\]
Let $u^{(0)} \in U_0 \left( C \left( X^{(0)}, \, M_n \right) \right)$
be a unitary such that
\[
{\ts{ \left\|
  \left[ u^{(0)} (x) \rsz{ \left[a (x)^* a (x) \right]^{1/2} }
            - a (x) \right] [1 - p (x)]  \right\| }} < \dt
\]
for $x \in X^{(0)}$.
Then there exists a unitary $u \in U_0 (C (X, M_n))$ such that
$u |_{X^{(0)}} = u^{(0)}$ and
\[
{\ts{ \left\| \left[ u (x) \rsz{ \left[a (x)^* a (x) \right]^{1/2} }
            - a (x) \right] [1 - p (x)]  \right\| }} < \ep
\]
for all $x \in X$.
Moreover, if we are given a homotopy $t \mapsto u_t^{(0)}$ from
$1$ to $u^{(0)}$ in $U \left( C \left( X^{(0)}, \, M_n \right) \right)$,
then $u$ can be chosen
such that there is a homotopy $t \mapsto u_t$ from
$1$ to $u$ in $U (C (X, M_n))$ such that $u_t |_{X^{(0)}} = u_t^{(0)}$.
\end{lem}

\begin{pff}
We may as well assume a homotopy $t \mapsto u_t^{(0)}$ is given.

Let $r$ be the least integer such that $r \geq \frac{1}{2} d$.

We choose $\dt$ by an inductive process.
Set $\dt_n = \ep$.
Given $\dt_{k + 1} > 0$, choose $\dt_k > 0$ so small that the value
$\dt = 2 \dt_k / \af$ works for $\ep = \frac{1}{2} \dt_{k + 1}$
in Lemma~\ref{D1}, and also so small that $\dt_k < \dt_{k + 1}$.
Then set $\dt = \dt_{r - 1}$.

Now we start the proof.
First observe that $\rank (p (x)) \geq r$ for all $x \in X$.
Indeed,
$\left\| \rsz{ \left[a (x)^* a (x) \right]^{1/2} } \xi \right\|
                   = \| a (x) \xi \|$
for all $x$ and $\xi \in \C^n$.
Therefore
$\left\| \rsz{ \left[a (x)^* a (x) \right]^{1/2} } \xi \right\|
                            < \af \|\xi\|$
for $\xi \in E_x \setminus \{ 0\}$.
Let $F_x$ be the linear span of the eigenspaces of
$\left[a (x)^* a (x) \right]^{1/2}$ for eigenvalues in $[ \af, \infty)$.
Then
$\left\| \rsz{ \left[a (x)^* a (x) \right]^{1/2} } \xi \right\|
                \geq \af \| \xi \|$
for all $\xi \in F_x$.
So $E_x \cap F_x = \{ 0 \}$.
It follows that $\rank (p (x)) = n - \dim (F_x) \geq \dim (E_x) \geq r$.

Applying Proposition~\ref{C9},
write $p (x) = \sup ( \{p_k (x) \colon  x \in U_k \})$,
where $p_k$ has rank $k$, the $U_k$ are open,
$\bigcup_{k = 0}^n U_k = X$, and
$p_k (x) \leq p_l (x)$ for $k \leq l$ and $x \in U_k \cap U_l$.
Since $\rank (p (x)) \geq r$ for all $x \in X$,
\wolog\  $U_k = \varnothing$ for $k < r$.
We now construct, by induction on $k$, closed sets $Y_k \subset X$
such that
\[
X^{(0)} \cup \{ x \in X \colon \rank (p (x)) \leq k \} \subset Y_k \subset
  X^{(0)} \cup \bigcup_{l = r}^k U_l
\]
and
\[
\{ x \in X \colon \rank (p (x)) \leq k \} \subset \sint (Y_k),
\]
unitaries $v_k \in U_0 (C (Y_k, \, M_n))$ such that
$v_k |_{X^{(0)}} = u^{(0)}$ and
\[
{\ts{ \left\|  \left[ v_k (x) \rsz{ \left[a (x)^* a (x) \right]^{1/2} }
            - a (x) \right] [1 - p (x)]  \right\| }} < \dt_k
\]
for all $x \in Y_k$,
and unitary homotopies
$(t, x) \mapsto w_t^{(k)} (x)$ in $C (Y_k, \, M_n)$ with
\[
w_0^{(k)} = 1, \,\,\,\,\,\,  w_1^{(k)} = v_k, \andeqn
w_t^{(k)} |_{X^{(0)}} = u_t^{(0)}.
\]

We start by taking $Y_{r - 1} = X^{(0)}$, $v_{r - 1} = u^{(0)}$,
and $w_0^{(r - 1)} = u_t^{(0)}$.

Suppose now we are given $Y_k$, $v_k$, and $w_0^{(k)}$.
Let
\[
R = \{ x \in X \colon \rank (p (x)) = k + 1\} \setminus \sint (Y_k).
\]
We need two facts about $R:$ that it is closed, and that the union of
$Y_k$ and any \nbhd\  of $R$ is a \nbhd\  of
$\{ x \in X \colon \rank (p (x)) \leq k + 1\}$.
For the first, let $(x_{\ld})$ be a net in $R$ with $x_{\ld} \to x$.
One easily sees that
\[
\rank (p (x)) \leq \liminf \rank (p (x_{\ld})) \leq k + 1.
\]
If now $\rank (p (x)) < k + 1$, then $x \in \sint (Y_k)$,
by the assumption on $Y_k$.
This is a contradiction.
So $\rank (p (x)) = k + 1$, and $x \in R$ because $\sint (Y_k)$ is open.
For the second, let $Z$ be a \nbhd\  of $R$.
Then $\sint (Z)$ contains $R$, and $\sint (Y_k)$ contains all other
points $x \in X$ such that $\rank (p (x)) = k + 1$.
Therefore $\sint (Y_k \cup Z)$ contains $R$.
The assumption on $Y_k$ implies that
$\sint (Y_k)$ contains $\{ x \in X \colon \rank (p (x)) \leq k \}$,
so that $\sint (Y_k \cup Z)$ contains
$\{ x \in X \colon \rank (p (x)) \leq k + 1\}$.

Clearly $U_{k + 1}$ is a \nbhd\  of $R$, and $p_{k + 1} (x) = p (x)$
for $x \in R$.
Therefore
\[
[1 - p_{k + 1} (x)] \left[a (x)^* a (x)\right]^{1/2} [1 - p_{k + 1} (x)]
  = [1 - p (x) ]  \left[a (x)^* a (x) \right]^{1/2} [1 - p (x)]
\]
is invertible in $[1 - p_{k + 1} (x) ] M_n [1 - p_{k + 1} (x) ]$
for $x$ in $R$.
So
\[
c_0 (x) = a (x)
      \left( [1 - p_{k + 1} (x)]  \left[a (x)^* a (x) \right]^{1/2}
                 [1 - p_{k + 1} (x)] \right)^{-1}
\]
(inverse taken in $[1 - p_{k + 1} (x) ] M_n [1 - p_{k + 1} (x) ]$)
exists for $x$ in some \nbhd\  $V$ of $R$.
Moreover,
$\| c_0 (x) [1 - p_{k + 1} (x)]\xi\| = \| [1 - p_{k + 1} (x)] \xi\|$
for $\xi \in \C^n$ and $x \in R$.
Therefore, using the compactness of the
closed unit ball of $\C^n$, we may reduce the size of $V$ so
that in addition the partial isometry
$c (x) = c_0 (x)  \left[c_0 (x)^* c_0 (x) \right]^{-1/2}$,
with initial \pj\  $1 - p_{k + 1} (x)$,
is defined for $x \in V$.
(Functional calculus is taken in
$[1 - p_{k + 1} (x) ] M_n [1 - p_{k + 1} (x) ]$.)

For $x \in Y_k$, we have
\[
{\ts{ \left\|  \left[ v_k (x) \rsz{ \left[a (x)^* a (x) \right]^{1/2} }
     - a (x) \right] [1 - p (x)]  \right\| }} < \dt_k,
\]
and for $x \in R$ we have
\[
{\ts{ \left[ c (x) \rsz{ \left[a (x)^* a (x) \right]^{1/2} }
     - a (x) \right] }} [1 - p (x)] = 0
\]
(using the relation $p (x) = p_{k + 1} (x)$ and the fact that $p (x)$
commutes with\linebreak
$\left[a (x)^* a (x) \right]^{1/2}$).
Using again $p (x) = p_{k + 1} (x)$ for $x \in R$, we get
\[
{\ts{ \left\| [ v_k (x) - c (x) ]
      \rsz{ \left[a (x)^* a (x) \right]^{1/2} }
              [1 - p_{k + 1} (x)] \right\| }} < \dt_k
\]
for $x \in Y_k \cap R$.
Now $\left[a (x)^* a (x) \right]^{1/2}$ commutes with
$1 - p_{k + 1} (x) = 1 - p (x)$, and
\[
{\ts{ \left\|  \rsz{ \left( [1 - p_{k + 1} (x)]
    \rsz{ \left[a (x)^* a (x) \right]^{1/2} }
         [1 - p_{k + 1} (x)] \right)^{-1} } \right\| }} <
  \frac{1}{\af}
\]
(inverse taken in $[1 - p_{k + 1} (x) ] M_n [1 - p_{k + 1} (x) ]$), so
\[
{\ts{ \left\| \rsz{ [ v_k (x) - c (x) ]
                      [1 - p_{k + 1} (x)] } \right\| }} <
         \frac{\dt_k}{\af}
\]
for $x \in Y_k \cap R$.
Choose a \nbhd\  $W$ of $Y_k \cap R$ with $W \subset V$ such that
this last norm is at most $2 \dt_k / \af$ for $x \in W$.
Then $(V \setminus Y_k) \cup W$ is a \nbhd\  of $R$.
Choose a closed \nbhd\  $Z$ of $R$ with
$Z \subset (V \setminus Y_k) \cup W$, and so small that
\[
{\ts{ \left\| \left[ c (x) \rsz{ \left[a (x)^* a (x) \right]^{1/2} }
     - a (x) \right] [1 - p_{k + 1} (x)] \right\| }} <
               {\textstyle{\frac{1}{2}}} \dt_{k + 1}
\]
for $x \in Z$.

Apply Lemma~\ref{D1} to $s_1 (x) = v_k (x) [1 - p_{k + 1} (x)]$ on
$Z_1 = Z \cap Y_k$ and $s_2 (x) = c (x)$ on $Z_2 = Z$.
The choice of $\dt_k$ provides a partial isometry $s$ with
initial \pj\  $(1 - p_{k + 1}) |_{Z}$ such that
\[
s |_{Z \cap Y_k} = [v_k (1 - p_{k + 1} )] |_{Z \cap Y_k} \andeqn
       \| s - c \| < {\textstyle{\frac{1}{2}}} \dt_{k + 1}.
\]
Let $q \in C (Z, M_n)$ be the \pj\  $q = 1 - s s^*$.

Apply Proposition~4.2~(2) of \cite{PhN}
with
$Z$ in place of $X$, with $Z \cap Y_k$ in place of $Y$, with
\[
p_1 = p_{k + 1}, \,\,\,\,\,\, p_2 = q, \,\,\,\,\,\,
q_1 = 1 - p_{k + 1}, \,\,\,\,\,\, q_2 = 1 - q, \,\,\,\,\,\,
v_0 = v_k p_{k + 1} |_{Z \cap Y_k}, \andeqn u = 1,
\]
with $s$ as given, and with $w_t^{(k)} |_{Z \cap Y_k}$ in place of
$w_t^{(0)}$.
This gives $v \in C (Z, M_n)$ with
\[
v^* v = p_{k + 1}, \,\,\,\,\,\, v v^* =  q, \andeqn
     v |_{Z \cap Y_k} = v_0,
\]
and a homotopy $(t, x) \mapsto w_t (x)$
of unitaries in $C (Z, M_n)$ such that
\[
w_0 = 1, \,\,\,\,\,\, w_1 = s + v, \andeqn
           w_t |_{Z \cap Y_k} = w_t^{(k)} |_{Z \cap Y_k}.
\]

Define $Y_{k + 1} = Z \cup Y_k$, and set
\[
v_{k + 1} (x) = \left\{ \begin{array}{lcl}
     v_k (x)             &  &  x \in Y_k  \\
     s (x) + v (x)       &  &  x \in Z.
         \end{array} \right.
\]
Since $Z$ contains a \nbhd\  of $R$, the facts about $R$ discussed right
after its choice imply that $Y_{k + 1}$ contains a \nbhd\  of
$\{ x \in X \colon \rank (p (x)) \leq k + 1 \}$.
By construction, we have
$Y_{k + 1} \subset  X^{(0)} \cup \bigcup_{l = r}^{k + 1} U_l$.
Define
\[
w_t^{(k + 1)} (x) = \left\{ \begin{array}{lcl}
     w_t^{(k)} (x)             &  &  x \in Y_k  \\
     w_t       (x)             &  &  x \in Z.
         \end{array} \right.
\]
Then $w_t^{(k + 1)}$ is a homotopy from $1$ to $v_{k + 1}$ in
$U ( C (Y_{k + 1}, \, M_n))$, such that\linebreak
$w_t^{(k + 1)} |_{Y_{k + 1}} = w_t^{(k)}$.
Moreover,
\[
v_{k + 1} |_{X^{(0)}} = v_k |_{X^{(0)}} = u^{(0)} \andeqn
w_t^{(k + 1)} |_{X^{(0)}} = w_t^{(k)} |_{X^{(0)}} = u_t^{(0)}.
\]

It remains only to show that
\[
{\ts{ \left\|  \left[ v_{k + 1} (x)
     \rsz{ \left[a (x)^* a (x) \right]^{1/2} }
            - a (x) \right] [1 - p (x)]  \right\| }} < \dt_{k + 1}
\]
for $x \in Y_{k + 1}$.
This estimate holds on $Y_k$ because $v_{k + 1} |_{Y_k} = v_k$
and $\dt_k < \dt_{k + 1}$.
For $x \in Z$, we have (using $p_{k + 1} (x) \leq p (x)$
and $v_{k + 1} |_Z = s$)
\begin{align*}
 & {\ts{ \left\|  \left[ v_{k + 1} (x)
       \rsz{ \left[a (x)^* a (x) \right]^{1/2} }
            - a (x) \right] [1 - p (x)]  \right\| }}  \\
 & \hspace{3em} \leq
 {\ts{ \left\|  \left[ s (x) \rsz{ \left[a (x)^* a (x) \right]^{1/2} }
            - a (x) \right] [1 - p_{k + 1} (x)]  \right\| }}  \\
 & \hspace{3em} \leq
\| s (x) - c (x) \| \| a \| +
 {\ts{ \left\|  \left[ c (x) \rsz{ \left[a (x)^* a (x) \right]^{1/2} }
            - a (x) \right] [1 - p_{k + 1} (x)]  \right\| }}  \\
 & \hspace{3em} < {\textstyle{\frac{1}{2}}} \dt_{k + 1}
                         + {\textstyle{\frac{1}{2}}} \dt_{k + 1}
              = \dt_{k + 1}.
\end{align*}
This completes the induction.

The proof is now finished by setting $u = v_n$ and $u_t = w_t^{(n)}$.
\end{pff}

\medskip

We can now prove the result on approximate polar decomposition.

\begin{prp}\label{D3}
Let $A$ be a \rsha\  with total space $X$, and let $m$ and $d$ be
its matrix size and \tdim\  functions as in Definition~\ref{A2}.
Let $\af, \, \ep > 0$.
Let $a \in A$, and suppose that for every $x \in X$ there is a subspace
$E_x$ of $\C^{m (x)}$
with $\dim (E_x) \geq \frac{1}{2} d (x)$ and  such that
$\| \ev_x (a) \xi \| < \af \| \xi \|$ for $\xi \in E_x \setminus \{ 0\}$.
Then there is a unitary $u \in U_0 (A)$ such that
$\left\| u \rsz{ \left(a^* a\right)^{1/2} } - a \right\| < 2 \af + \ep$.
\end{prp}

\begin{pff}
We first to reduce to the case in which,
following the notation of Definition~\ref{A2},
the first space $X_0$ has just one point.
To do this, replace $C \left( X_0, \, M_{n (0)} \right)$ by
$M_{n (0)} \oplus_{M_{n (0)}} C \left( X_0, \, M_{n (0)} \right)$,
where the map
$M_{n (0)} \to M_{n (0)}$ is the identity map and the map
$C \left( X_0, \, M_{n (0)} \right) \to M_{n (0)}$ is $\ev_x$
for some $x \in X_0$.
This change increases the length of the decomposition by $1$,
but does not affect any of the hypotheses, or the conclusion, of
the proposition.

Assuming now that all decompositions we consider start with a
one point space, we prove the following by induction on the length:
Let $a \in A$ satisfy the hypotheses, and let
$p (x) = \ch_{ ( -\infty, \af)}
   \left( \ev_x \left(  (a^* a )^{1/2} \right) \right)$
for $x \in X$.
Then there is a unitary $u \in U_0 (A)$ such that
\[
{\ts{ \left\| \ev_x \left( u \rsz{ \left( a^* a\right)^{1/2} }
                           - a \right)
                      [1 - p (x)]  \right\| }}
   < \ep
\]
for all $x \in X$.
This unitary will then be shown to satisfy
the conclusion of the proposition.

The case of length zero is now trivial, since then $A = M_{n (0)}$.
So assume that the result is known for length $L$.
By scaling $a$, $\af$, and $\ep$, \wolog\  $\| a \| \leq 1$.
Let $A = B \oplus_{C \left( X^{(0)}, \, M_n \right)} C (X, M_n)$, with
$\ph \colon B \to C \left( X^{(0)}, \, M_n \right)$ unital and
$\rh \colon C (X, M_n) \to C \left( X^{(0)}, \, M_n \right)$
the restriction map, where $B$ is a \rsha\  of length $L$.
Let $Y$ be the total space of $B$, so that
the  total space of $A$ is $Y \coprod X$.
Let $d = \dim (X)$.
Let $b$ be the image of $a$ in $B$, and let $a_0$ be the
image of $a$ in $C (X, M_n)$, under the obvious \pj\  maps.
Note that
$\ch_{ ( -\infty, \af)}
   \left( \ev_y \left(  (b^* b )^{1/2} \right) \right) = p (y)$
for  $y \in Y$.
For the given values of $\af, \, \ep$, and $n$, let $\dt$ be as in
Lemma~\ref{D2}.
Choose, using the induction assumption,  a unitary $v \in B$ such that
\[
{\ts{ \left\| \ev_y \left( v \rsz{ \left( b^* b \right)^{1/2} }
                                        - b \right)
                        [1 - p (y)]  \right\| }}
   < \dt
\]
for all $y \in Y$, and such that there is a homotopy
$t \mapsto v^{(t)}$ in $U (B)$ with $v^{(0)} = 1$ and $v^{(1)} = v$.

For each $x \in X^{(0)}$, the map $c \to \ph (c) (x)$ is a
finite dimensional representation of $B$.
It follows from Lemma~2.1 of \cite{PhN}
that there are  a unitary $w \in M_n$
and points $y_1, \dots, y_l \in Y$ such that
\[
\ph (c) (x) = w \left( \bigoplus_{j = 1}^l \ev_{y_j} (c) \right) w^*
\]
for all $c \in B$.
So
\[
a_0 (x) = \ev_x (a)
            = w \left( \bigoplus_{j = 1}^l \ev_{y_j} (b) \right) w^*
           \andeqn
\ph (v) (x) = w \left( \bigoplus_{j = 1}^l \ev_{y_j} (v) \right) w^*.
\]
Applying functional calculus to the first of these
(with $(a^* a )^{1/2}$ and $(b^* b )^{1/2}$ in place of $a$ and $b$)
gives
\[
p (x) = w \left( \bigoplus_{j = 1}^l p (y_j) \right) w^*.
\]
It now follows from the choice of $v$ that
\[
{\ts{ \left\| \left[ \ph (v) (x)
     \rsz{ \left[ a_0 (x)^* a_0 (x) \right]^{1/2} }
                              - a_0 (x) \right] [1 - p (x)]  \right\| }}
              < \dt
\]
for all $x \in X^{(0)}$.
Moreover, $t \mapsto \ph \left( v^{(t)} \right)$ is a homotopy
from $1$ to $\ph (v)$.
According to Lemma~\ref{D2} and the choice of $\dt$,
there exist a unitary
$u_0 \in C (X, M_n)$ and a unitary homotopy $t \mapsto u_0^{(t)}$
from $1$ to $u_0$, such that $u_0 |_{X^{(0)}} = \ph (v)$ and
\[
{\ts{ \left\| \left[ u (x) \rsz{ \left[a_0 (x)^* a_0 (x) \right]^{1/2} }
            - a_0 (x) \right] [1 - p (x)]  \right\| }} < \ep
\]
for all $x \in X$, and such that
$u_0^{(t)} |_{X^{(0)}} = \ph \left( v^{(t)} \right)$.
Then $u = (v, u_0)$ is a unitary in $U_0 (A)$ such that
\[
{\ts{
\left\| \ev_x \left( u \rsz{ \left( a^* a\right)^{1/2} } - a \right)
                                    [1 - p (x)]  \right\| }}
   < \ep
\]
for all $x \in Y \coprod X$.
Moreover, the homotopy
$t \mapsto u^{(t)} = \left( v^{(t)}, \rsz{ u_0^{(t)} }  \right)$ in
$U (A)$ shows that $u \in U_0 (A)$.
This completes the induction, and the proof of the claim.

To get the desired estimate, write
\begin{align*}
& {\ts{ \left\| \ev_x \left( u \rsz{ \left( a^* a\right)^{1/2} }
        - a \right) \right\| }} \\
& \hspace{3em} \leq
 {\ts{ \left\| \ev_x \left( u \rsz{ \left( a^* a\right)^{1/2} }
        - a \right) [1 - p (x)]  \right\| }}
    + {\ts{ \left\| \ev_x \left( u \rsz{ \left( a^* a\right)^{1/2} }
                   - a \right) p (x)  \right\| }}.
\end{align*}
The first term on the right is less than $\ep$.
For the second, we have, for all $\xi$,
\[
{\ts{ \left\| \ev_x \left( u \rsz{ \left( a^* a\right)^{1/2} } \right)
                    p (x) \xi \right\| }}
 = {\ts{ \left\| \ev_x \left( \rsz{ \left( a^* a\right)^{1/2} } \right)
                    p (x) \xi \right\| }}
  = \left\| \ev_x (a) p (x) \xi \right\| \leq \af \| \xi \|,
\]
so
\[
{\ts{ \left\| \ev_x \left(
    u \rsz{ \left( a^* a \right)^{1/2} } - a \right) p (x)  \right\| }}
  \leq 2 \af \| \xi \|.
\]
It follows that
$\left\| \ev_x \left( u \rsz{ \left( a^* a\right)^{1/2} }
                             - a \right) \right\|
                      <  \ep + 2 \af$.
This is true for all $x$ in the total space of $A$, so
$\left\| u \rsz{ \left( a^* a\right)^{1/2} } - a \right\| <  \ep + 2 \af$.
\end{pff}

\medskip

In order to apply our approximate polar decomposition to simple
direct limits, we need the following lemma.

\begin{lem}\label{D4}
Let $A = \dirlim A_i$ be an infinite dimensional
simple direct limit of \rsha s, such that
all the maps $\ph_{i j} \colon A_i \to A_j$ in the system are unital and
injective.
Let $X_i$ be the total space of $A_i$, and let $m_i \colon X_i \to \N$ be
the matrix size function.
Let $a \in A_i$, for some $i$, be noninvertible.
Then for every $n \in \N$ and $\ep > 0$ there exists $j_0$ such that,
for every $j \geq j_0$ and every
$x \in X_j$ there is a subspace $E_x$ of $\C^{m_j (x)}$
such that $\dim (E_x) \geq n$ and
$\| \ev_x ( \ph_{i j} (a)) \xi \| < \ep \| \xi \|$
for $\xi \in E_x \setminus \{ 0\}$.
\end{lem}

\begin{pff}
Since
$\left\| \rsz{ \left(b^* b \right)^{1/2} } \xi \right\| = \| b \xi \|$
for any $b$ and $\xi$, we may replace $a$ by $\left(a^* a\right)^{1/2}$,
and thus assume $a \geq 0$.

Let $f \colon [0, \infty) \to [0, 1]$ be a \ct\  function such that
$f (0) \neq 0$ but $f (t) = 0$ for $t \geq \frac{1}{2} \ep$.
Then $f (a)$ is a nonzero selfadjoint element of $A_i$.
Applying Lemma~\ref{G8},
we obtain $j_0$ such that, for every $j \geq j_0$ and
every  $x \in X_j$, the matrix $\ev_x ( \ph_{i j} ( f (a)))$ has
rank at least $n$.
Let $g \colon [0, \infty) \to [0, 1]$ be a \ct\  function such that
$g = 1$ on $\left[ 0, \frac{1}{2} \ep \right]$ and $g = 0$ on $[ \ep, \infty)$.
Since $\| a g (a) \| < \ep$ and $g (a) f (a) = f (a)$, it
follows that $\| \ev_x ( \ph_{i j} (a)) \xi \| < \ep \| \xi \|$
for $\xi \in \ev_x ( \ph_{i j} ( f (a))) \C^{n_j (x)} \setminus \{ 0\}$.
The lemma is therefore proved by taking
$E_x = \ev_x ( \ph_{i j} ( f (a))) \C^{m_j (x)}$.
\end{pff}

\begin{thm}\label{D5}
Let $A = \dirlim A_i$ be a simple direct limit of separable \rsha s.
Assume the system has no dimension growth, that is, there is
$d \in \N$ such that all $A_i$ have \tdim\  at most $d$.
Then $\tsr (A) = 1$.
\end{thm}

\begin{pff}
If $A$ is finite dimensional, the conclusion is obvious.
So assume $A$ is infinite dimensional.

We first consider the case in which
all the maps $\ph_{i j} \colon A_i \to A_j$ in the system are unital and
injective.
By Lemma~3.5 of \cite{DNNP}, it suffices to let $a \in A_i$ for some
$i$, let $\ep > 0$, and find $j \geq i$ and an invertible element
$c \in A_j$ such that $\| \ph_{i j} (a) - c \| < \ep$.
Use Lemma~\ref{D4} to find $j \geq i$ such that for every $x \in X_j$
there is a subspace $E_x$ of $\C^m$ (where $m$ is the
matrix size of $A_j$ at the point $x$) such that
$\dim (E_x) \geq \frac{1}{2} d$
and $\| \ev_x ( \ph_{i j} (a)) \xi \| < \frac{1}{4} \ep \| \xi \|$
for $\xi \in E_x \setminus \{ 0\}$.
Set $b = \ph_{i j} (a)$.
Using Proposition~\ref{D3}, find a unitary $u \in A_j$ such that
$\left\| u \rsz{ \left(b^* b \right)^{1/2} } - b \right\|
           < \frac{3}{4} \ep$.
Then
$c = u \left[ \rsz{ \left(b^* b \right)^{1/2} }
                     + \frac{1}{4} \ep \cdot 1 \right]$
is an invertible element of $A_j$ which satisfies
$\| \ph_{i j} (a) - c \| = \| b - c \| < \ep$.

Now we drop the assumption that the maps are unital.
\Wolog\  $A_0 \neq 0$.
Let $p_i = \ph_{0 i} (1_{A_0})$, which is a nonzero \pj\  in $A_i$,
and let $p$ be the image of $1_{A_0}$ in $A$.
Then $p A p \cong \dirlim p_i A_i p_i$.
By Corollary~1.11 of \cite{PhN},
the algebras $p_i A_i p_i$ are \rsha s
with \tdim\  at most $d$.
Also $p A p$ is simple because $A$ is.
Therefore the previous case shows $\tsr (p A p) = 1$.
Since  $p A p$ is stably isomorphic to $A$ (Theorem~2.8 of \cite{Br}),
it follows from Theorem~3.6 of \cite{Rf1} that $\tsr (A) = 1$.

Finally, we use Proposition~\ref{G10} to
drop the assumption that the maps are injective.
\end{pff}

\section{Examples and applications}

In this section, we apply our general results to crossed products by
minimal homeomorphisms.
In particular, we compute the Elliott invariant \cite{El5} of
crossed products by \mh s of odd spheres of dimension at least $3$,
and we give a much more direct calculation of the
Elliott invariant of the crossed product by a Furstenberg transformation
on the $2$-torus (see \cite{Kd}).
We also give a related example which shows the
failure of the generalization to direct limits of \rsha s
of results on Riesz decomposition and real rank zero.
The applications and the example use the subalgebras of the crossed
product studied in \cite{Ln} (see Example~1.6 of \cite{PhN}).

We note here that the preprint \cite{Ln} is expected to be absorbed
into \cite{LP2}, in which a much stronger result is proved.
(However, a sketch of what is needed here has been published in
Section~3 of \cite{LP}.)
For our purposes, the main result of \cite{Ln} is that if $h$ is a
\mh\  of a finite dimensional infinite compact metric space $X$,
then a useful ``large'' subalgebra
($A_{\{ x \} }$, described below) of the crossed product
$C^* (\Z, X, h)$ is a simple direct limit of \rsha s
with no dimension growth.
In case $X$ is a \mf\  and $h$ is actually a diffeomorphism,
it is shown in \cite{LP2} that $C^* (\Z, X, h)$ is itself
such a direct limit.
The application of the theorems of this paper to $C^* (\Z, X, h)$
gives considerably more information than their application to the
subalgebra $A_{\{ x \} } \subset C^* (\Z, X, h)$.
We have two reasons for giving the weaker results here.
First, it is quite straightforward to prove that $A_{\{ x \} }$
is a simple direct limit of \rsha s with no dimension growth,
while the proof of the corresponding theorem for the crossed product
is extremely long.
The result for $A_{\{ x \} }$ already gives enough information to
compute the Elliott invariant of the crossed product.
Second, the direct limit decomposition for the crossed product has only
been proved for minimal diffeomorphisms; if $X$ is not a manifold,
or if it is but $h$ is not smooth, the results presented
here are the best currently known.

We begin by recalling from \cite{Ln} the
subalgebra $A_{\{ x \} }$ of the crossed product, and stating its
relation with the entire crossed product algebra.

\begin{thm}\label{Q1}
Let $X$ be an infinite compact metric space,
let $h$ be a minimal homeomorphism of $X$, and set $A = C^* (\Z, X, h)$.
Let $u \in A$ be the unitary representing the generator of $\Z$.
For $x \in X$, set
\[
A_{\{x\}} = C^* ( C (X), \, u C_0 (X \setminus \{x\}))
                  \subset C^* (\Z, X, h),
\]
as in \cite{Ln} (also see Example~1.6 of \cite{PhN}),
and let $\io \colon A_{\{x\}} \to A$ be the inclusion map.
Then:
\bit
\item[(1)]
$A_{\{x\}}$ is simple.
\item[(2)]
$A_{\{x\}}$ is a direct limit of a system of \rsha s with
\tdim\  at most $\dim (X)$, in which the
maps are all unital and injective.
\item[(3)]
$\io_* \colon K_0 \left( A_{\{x\}} \right) \to K_0 ( A)$ is an isomorphism.
\item[(4)]
There is a short exact sequence
\[
0 \longrightarrow K_1 {\ts{ \left( A_{\{x\}} \right) }}
            \stackrel{\io_*}{\longrightarrow}
            K_1 ( A) \stackrel{\gm}{\longrightarrow} \Z
            \longrightarrow 0,
\]
in which $\gm ([u]) = 1$.
\item[(5)]
There is a one to one correspondence between normalized traces $\ta$
on $A$ and $h$-invariant Borel probability measures $\mu$ on $X$,
given by
\[
\ta \left(\sum_{k} f_k u^k \right) = \int_X f_0 \, d \mu.
\]
Moreover, the map $\ta \mapsto \ta \circ \io$ is a bijection between
normalized traces $\ta$ on $A$ and normalized traces on $A_{\{x\}}$.
\eit
\end{thm}

\begin{pff}
(1) This is Proposition~12 of \cite{Ln}.

(2) This is Example~1.6 of \cite{PhN}
(derived from Theorem~3 of \cite{Ln}).

(3), (4) Example~2.6 of \cite{Pt} gives an exact sequence
\[
\begin{array}{ccccc}
K_0 (K (l^2 (\Z)))  & \longrightarrow &
    K_0 {\ts{ \left( A_{\{x\}} \right) }}  &
    \longrightarrow & K_0 (A)   \\
\uparrow  & & & & \downarrow       \\
K_1 (A) & \longleftarrow & K_1 {\ts{ \left( A_{\{x\}} \right) }} &
      \longleftarrow & K_1 (K (l^2 (\Z)))
\end{array}
\]
in which the maps $K_* {\ts{ \left( A_{\{x\}} \right) }}  \to K_* (A)$
are induced by the inclusions.
Moreover, we claim that the map $K_1 (A) \to K_0 (K (l^2 (\Z)))$
sends the class $[u]$ of
the generating unitary to a generator of $K_0 (K (l^2 (\Z)))$.
Since $K_0 (K (l^2 (\Z))) \cong \Z$ and $K_1 (K (l^2 (\Z))) = 0$,
the conclusions will then follow from exactness.

To prove the claim, we chase through the definitions in \cite{Pt}.
(The map in question is called $[L]_*$ there, and it is defined in
the discussion following Lemma~3.10 of \cite{Pt}.)
We find that this map is determined by the odd Kasparov
$A$--$K (l^2 (\Z))$-bimodule
$\left( K (l^2 (\Z)), \, \overline{\ld}, \, 1 - 2 p \right)$,
in which the right $K (l^2 (\Z))$ action on $K (l^2 (\Z))$ is simply
multiplication, in which
$\overline{\ld} \colon A \to L (l^2 (\Z)) = M ( K (l^2 (\Z)))$ is the
regular representation induced as in 7.7.1 of \cite{Pd0}
by the representation $\ev_x$ of $C (X)$, and in which $p$ is the
\pj\  from $l^2 (\Z)$ onto $l^2 (\{ 1, 2, \dots \})$.
Since $\overline{\ld} (u)$ is the bilateral shift, this map does
indeed send $[u]$ to a generator of $K_0 (K (l^2 (\Z)))$.

(5) This is Proposition~16 of \cite{Ln} and the well known
correspondence between traces on $A$ and $h$-invariant
Borel probability measures on $X$.
\end{pff}

\begin{thm}\label{Q2}
Let $X$, $h$, and $A_{\{x\}}$ be as in Theorem~\ref{Q1},
and assume in addition that $\dim (X) < \infty$.
Then:
\bit
\item[(1)]
$\tsr {\ts{ \left( A_{\{x\}} \right) }} = 1$.
\item[(2)]
The map
$U {\ts{ \left( A_{\{x\}} \right) }} /
              U_0 {\ts{ \left( A_{\{x\}} \right) }}
    \to K_1 {\ts{ \left( A_{\{x\}} \right) }}$
is an isomorphism.
\item[(3)]
The \pj s in $\Mi {\ts{ \left( A_{\{x\}} \right) }}$ satisfy cancellation.
\item[(4)]
The algebra $A_{\{x\}}$ satisfies
Blackadar's Second Fundamental Comparability Question:
if $p, \, q \in \Mi {\ts{ \left( A_{\{x\}} \right) }}$ are \pj s such that
$\ta (p) < \ta (q)$ for all normalized traces
$\ta$ on $A_{\{x\}}$, then $p \precsim q$.
\item[(5)]
$K_0 {\ts{ \left( A_{\{x\}} \right) }}$
is unperforated for the strict order:
if $n \et > 0$ in $K_0 {\ts{ \left( A_{\{x\}} \right) }}$, with $n > 0$,
then $\et > 0$.
\eit
\end{thm}

\begin{pff}
Parts~(1) and~(2) of Theorem~\ref{Q1} show that $A_{\{x\}}$ is a simple
direct limit, with
no dimension growth, of a system of separable \rsha s.
Therefore part~(1) follows from Theorem~\ref{D5},
part~(2) follows from Corollary~\ref{G9} and Theorem~\ref{C1},
part~(3) follows from Corollary~\ref{G9} and Theorem~\ref{C2},
part~(4) follows from Corollary~\ref{G9} and Theorem~\ref{C3},
and part~(5) follows from Corollary~\ref{G9} and Theorem~\ref{C4}.
\end{pff}

\medskip

While not all of the properties in this theorem can be transferred
to the entire crossed product $A$,
Theorem~\ref{Q1} does give us some information.

\begin{thm}\label{Q3}
Let $X$ be a finite dimensional infinite compact metric space,
and let $h$ be a minimal homeomorphism of $X$.
Then the map
\[
U (C^* (\Z, X, h)) / U_0 (C^* (\Z, X, h))
                 \longrightarrow K_1 (C^* (\Z, X, h))
\]
is surjective.
\end{thm}

\begin{pff}
Theorem~\ref{Q1}~(4) implies (using the notation there)
that $K_1 (C^* (\Z, X, h))$ is generated
by $\io_* \left( K_1 \left( A_{\{x\}} \right) \right)$ and $[u]$.
The image of the map
$U (C^* (\Z, X, h)) \to K_1 (C^* (\Z, X, h))$
is a subgroup of $K_1 (C^* (\Z, X, h))$ which contains
$\io_* \left( K_1 \left( A_{\{x\}} \right) \right)$
by Theorem~\ref{Q2}~(2),
and which obviously contains $[u]$,
so it is all of $K_1 (C^* (\Z, X, h))$.
\end{pff}

\begin{thm}\label{Q4}
Let $X$, $h$, $A$, and $A_{\{x\}}$ be as in Theorem~\ref{Q1},
and assume in addition that $\dim (X) < \infty$.
Then $\io_* \colon K_0 {\ts{ \left( A_{\{x\}} \right) }} \to K_0 ( A)$
is an {\emph{order}} isomorphism.
\end{thm}

\begin{pff}
By Theorem~\ref{Q1}~(3), we need only prove that $(\io_*)^{-1}$
is order preserving.
So let $\et \in K_0 (A)$ satisfy $\et > 0$.
Write $(\io_*)^{-1} (\et) = [p] - [q]$ for some \pj s
$p, \, q \in \Mi {\ts{ \left( A_{\{x\}} \right) }}$.
Let $\sm$ be any normalized trace on $A_{\{x\}}$.
By Theorem~\ref{Q1}~(5),
there is a normalized trace $\ta$ on $A$ such that
$\ta \circ \io = \sm$.
We have $\ta_* (\et) > 0$ because $A$ is simple and $\et > 0$.
Therefore $\ta ( \io (p)) > \ta (\io (q))$,
whence $\sm (p) > \sm (q)$.
Since $\sm$ is arbitrary,
Theorem~\ref{Q2}~(4) implies that $p \succsim q$, so that
$(\io_*)^{-1} (\et) = [p] - [q] > 0$.
\end{pff}

\begin{thm}\label{Q5}
Let $X$ be a finite dimensional infinite compact metric space,
and let $h$ be a minimal homeomorphism of $X$.
Then:
\bit
\item[(1)]
$C^* (\Z, X, h)$ satisfies the following K-theoretic version of
Blackadar's Second Fundamental Comparability Question:
if $\et \in K_0 (A)$ satisfies $\ta_* (\et) > 0$
for all normalized traces $\ta$ on $A$, then there is a
\pj\   $p \in \Mi (A)$ such that $\et = [p]$.
\item[(2)]
$K_0 (C^* (\Z, X, h))$ is unperforated for the strict order.
\eit
\end{thm}

\begin{pff}
This follows from Theorem~\ref{Q4}
and parts~(4) and~(5) of Theorem~\ref{Q2}.
\end{pff}

\medskip

We note at this point that the examples constructed by
Villadsen in \cite{Vl} do not satisfy~(2) (see Proposition~11~(ii) of
\cite{Vl}), and so also don't satisfy~(1).
Thus, the \ca s of minimal homeomorphisms of finite
dimensional compact metric spaces are not as badly behaved as
Villadsen's examples.
(In \cite{LP2}, we will show that crossed products by minimal
diffeomorphisms in fact have stable rank one.)
On the other hand, the \ca s covered by the real rank one
classification theorem of \cite{EGL2} have Riesz decomposition in $K_0$
even in the real rank one case (see Theorems~2.6 and 4.8 of \cite{LH}),
and we show in the next example that
the \ca s of minimal homeomorphisms of finite dimensional compact
metric spaces need not have Riesz decomposition in $K_0$.

\begin{exa}\label{Q8}
We compute the ordered $K_0$-group of Connes' example
(Example~4 in Section~5 of \cite{Cn})
of a simple unital stably finite \ca\  with no nontrivial \pj s.
As there, let $h$ be a minimal \diff\  of $S^3$.
Then $A = C^* (\Z, S^3, h)$ is simple,
and Corollary~3 in Section~5 of \cite{Cn}
implies that $\ta_* (K_0 (A)) = \Z$ for any normalized trace
$\ta$ on $A$.
(Also see Corollary~10.10.6 of \cite{Bl2}.)
Since $h$ has no fixed points, the Lefschetz Fixed Point Theorem
(Theorem~4.7.7 of \cite{Sp}) implies that $h$ is orientation preserving.
The Pimsner-Voiculescu exact sequence \cite{PV}
then implies that $K_0 (A) \cong \Z^2$.

We next note that all normalized traces on $A$ agree on $K_0 (A)$.
Indeed, if $\ta_1$ and $\ta_2$ are normalized traces, then
$t \mapsto [ (1 - t) \ta_1 + t \ta_2 ]_*$ is a homotopy, in an obvious
sense, of maps from $K_0 (A)$ to $\R$ whose ranges are contained in
$\Z$.
Therefore it is constant.
We let $\ta_*$ denote the map from $K_0 (A)$ to $\Z$ determined by
any trace.

The map $n \mapsto n \cdot [1]$, from $\Z$ to $K_0 (A)$, is a
left inverse of $\ta_*$.
It induces an isomorphism $K_0 (A) \cong \ker (\ta_*) \oplus \Z$,
with $\ker (\ta_*) \cong \Z$ also.
We may thus identify $K_0 (A)$ with $\Z^2$ in such a way that
$[1] = (0, 1)$ and $\ta_* (m, n) = n$.
Using Theorem~\ref{Q4}~(1),
we find that $(m, n) \geq 0$ exactly when $n > 0$ or
$(m, n) = (0, 0)$.

We note that $K_0 (A)$ does not have Riesz decomposition.
Indeed, write $(0, 2) = (1, 1) + (-1, 1)$, and note that
$(0, 1) \leq (0, 2)$
but there do not exist $\et, \, \mu \in K_0 (A)$ with
\[
(0, 0) \leq \et \leq (1, 1), \,\,\,\,\,\,
(0, 0) \leq \mu \leq (-1, 1), \andeqn \et + \mu = (0, 1).
\]
\end{exa}

\medskip

In this example, $S^3$ could be replaced by $S^n$ for any odd
$n \geq 3$.
Moreover, the outcome shows that the Elliott invariant depends only
on the space of $h$-invariant Borel probability measures.
The following conjecture is therefore a special case of the
Elliott classification conjecture.

\begin{cnj}\label{Q8.5}
Let $g$ and $h$ be minimal diffeomorphisms of $S^m$ and $S^n$
respectively, with $m, \, n \geq 3$ and odd.
Suppose the spaces of $g$-invariant and $h$-invariant
Borel probability measures are affinely homeomorphic.
Then there is an isomorphism of transformation group \ca s
$C^* (\Z, S^m, g) \cong C^* (\Z, S^n, h)$.
\end{cnj}

See Section~5 of \cite{LP} for further discussion of what the
Elliott classification conjecture might imply for crossed
products by \mh s, and in particular the contrast between
the suggested behavior of minimal diffeomorphisms of
high dimensional manifolds with what is known to happen
for \mh s of the Cantor set \cite{GPS} and of the circle.

The Connes example can also be used to show that some of the results
of \cite{BDR} and \cite{Gd} do {\emph{not}} generalize to direct limits
of \rsha s.
Examples of this type have been given previously.
See \cite{Th3} and \cite{JS} for the failure of real rank zero,
and \cite{Th2} (recalling from Theorem~10.17 and Proposition~2.1
of \cite{Gd0} that Riesz decomposition implies that the state space
of the group is a Choquet simplex)
for the failure of Riesz decomposition.
(More general results on existence of algebras with prespecified
Elliott invariants can be found in \cite{El2}.)
However, our example comes up naturally and has a simple proof.

\begin{exa}\label{Q9}
We give an example of a simple direct limit $B$ of separable \rsha s,
with no dimension growth, which has the following properties:
\bit
\item[(1)]
The \pj s distinguish the traces but $B$ does not
have real rank zero.
\item[(2)]
$K_0 (B)$ does not have Riesz decomposition.
\eit
Thus, there is no analog of Theorem~2 of \cite{BDR}, or of
Theorem~2.7 of \cite{Gd}, for direct limits of \rsha s.

Let $A$ be as in Example~\ref{Q8}, using
a uniquely ergodic minimal \diff\  $h$ of $S^3$.
Such a thing exists by \cite{FH}.
Then $A = C^* (\Z, S^3, h)$ is simple and has a unique normalized
trace $\ta$.
Let $x \in S^3$, and let $B = A_{\{x\}}$ be as in Theorem~\ref{Q1}.
As there, it is a simple unital direct limit of separable \rsha s,
with no dimension growth.
Moreover, by Theorem~\ref{Q1}~(5), it has a unique trace, since
there is a unique $h$-invariant Borel probability measure.
In particular, the \pj s distinguish the traces.
Corollary~3 in Section~5 of \cite{Cn}
implies that $A$ has no nontrivial \pj s.
Therefore $A_{\{x\}}$, being a subalgebra of $A$, also
has no nontrivial \pj s.
Consequently it does not have real rank zero.
The ordered $K_0$-group of $A_{\{x\}}$ is the same as for $A$, and
was shown in the previous example not to have Riesz decomposition.
\end{exa}

\medskip

As another example,
we compute the Elliott invariants of crossed products by
Furstenberg transformations on the $2$-torus,
recovering in particular the main result of \cite{Kd}
(namely, the order on $K_0$) with much less effort.
The computations of the (unordered) K-theory and
the effect of the traces on $K_0$ were done in \cite{Ji},
which was never published.
Moreover, there is now better machinery available \cite{Ex} for the
computation of the effect of the traces on $K_0$.
It therefore seems appropriate to start from scratch.

\begin{exa}\label{FT}
Fix $\te \in [0, 1] \setminus \Q$, a \ct\  function
$f_0 \colon S^1 \to \R$, and $n \in \Z \setminus \{ 0 \}$.
We define $h \colon S^1 \times S^1 \to S^1 \times S^1$ to be the inverse of
the homeomorphism
\[
( \zt_1, \zt_2 ) \mapsto
 {\ts{ \left( \rsz{ e^{2 \pi i \te} \zt_1, \,
     e^{2 \pi i f_0 (\zt_1)} \zt_1^n \zt_2 } \right) }}.
\]
(One sees that the given map does in fact have a
\ct\  inverse by writing down an explicit formula for it.
This homeomorphism is called $\ph_{f_0, \te}$ in \cite{Kd}.)
The homeomorphism $h$ is minimal by the remark after
Theorem~2.1 in Section~2.3 of \cite{Fr2}.
Normalized Lebesgue measure on $S^1 \times S^1$ is invariant,
and when $f_0$ is Lipschitz this is the only invariant measure
(Theorem~2.1 of \cite{Fr2}).
Define $\af \colon C ( S^1 \times S^1 ) \to C ( S^1 \times S^1 )$ by
$\af (f) = f \circ h^{-1}$.
Let
\[
A = C^* (\Z, \, S^1 \times S^1, \, h)
           = C^* (\Z, \, C ( S^1 \times S^1 ), \, \af).
\]

We compute the Elliott invariant of $A$, and we start by describing
the ingredients of the description we intend to prove.
Define $\af_0 \colon C (S^1) \to C (S^1)$ by
$\af_0 (f) (\zt) = f \left( e^{2 \pi i \te} \zt \right)$.
Then the crossed product $C^* (\Z, \, C (S^1), \, \af_0)$ is just the
irrational rotation algebra $A_{\te}$.
Moreover, the \hm\  $a \mapsto a \otimes 1$,
from
$C (S^1)$ to $C ( S^1) \otimes C (S^1 ) \cong C ( S^1 \times S^1 )$,
intertwines $\af_0$ and $\af$,
thus giving a \hm\  $\ph \colon A_{\te} \to A$.
Let $p \in A_{\te}$ be a \pj\  for which the unique trace $\ta$
on $A_{\te}$ satisfies $\ta (p) = \te$ (\cite{Rf0}).
Let $\io \colon C ( S^1 \times S^1 ) \to A$ be the inclusion,
and let $\bt \in K_0 (C ( S^1 \times S^1 ))$ be the Bott element.
We prove that the Elliott invariant of $A$ is given as follows:
\[
K_0 (A) \cong \Z \cdot [1] \oplus \Z \cdot \io_* (\bt)
   \oplus \Z \cdot [\ph (p)],
\]
every trace $\ta$ on $A$ satisfies
\[
\ta_* (m_1 [1] + m_2 \io_* (\bt) + m_3 [\ph (p)] ) = m_1 + m_3 \te,
\]
\[
K_0 (A)_+ = \{ m_1 [1] + m_2 \io_* (\bt) + m_3 [\ph (p)] \colon
  m_1 + m_3 \te > 0 \,\, {\mbox{or}} \,\, m_1 = m_2 = m_3 = 0 \},
\]
and
\[
K_1 (A) \cong \Z^3 \oplus \Z / n \Z.
\]

We first compute $K_* (A)$ using the Pimsner-Voiculescu exact sequence
\cite{PV}, which here takes the form
\[
\begin{array}{ccccc}
K_0 (C ( S^1 \times S^1 ))
    & \stackrel{\id - \af^{-1}_*}{\longrightarrow}
    & K_0 ( C ( S^1 \times S^1 ) )  & \longrightarrow & K_0 (A)   \\
{\mathrm{exp}} \uparrow \hspace*{2em}  & & & &
            \hspace*{1em} \downarrow \partial       \\
K_1 (A) & \longleftarrow & K_1 ( C ( S^1 \times S^1 ) ) &
    \stackrel{\id - \af^{-1}_*}{\longleftarrow} &
     K_1 (C ( S^1 \times S^1 ))
\end{array}.
\]
Let $z \in U ( C ( S^1))$ be $z (\zt) = \zt$.
Then we identify $K_1 ( C ( S^1 \times S^1 ) ) \cong \Z^2$ as the
free abelian group on the generators $[z \otimes 1]$ and $[1 \otimes z]$
(in that order).
We also identify $K_0 ( C ( S^1 \times S^1 ) ) \cong \Z^2$ as the
free abelian group on the generators $[1]$ and the Bott element $\bt$,
which we take to be the image of $[z] \otimes [z]$ under the
\hm\  %
\[
K_1 ( C (S^1)) \otimes K_1 ( C (S^1))
    \to K_0 ( C ( S^1 ) \otimes C (S^1 ) ).
\]
(See the beginning of Section~2 of \cite{Sc}.)
Again, for the purpose of writing group \hm s as matrices,
we take these generators in the order given.

To compute $\af_*$, we use the homotopic map given by
$f \mapsto f \circ h_0^{-1}$ with
\[
h_0^{-1} ( \zt_1, \zt_2 ) = ( \zt_1, \zt_1^n \zt_2 ).
\]
It is then clear that
\[
\af_* ([1]) = [1], \,\,\,\,\,\, \af_* ([z \otimes 1]) = [z \otimes 1],
\andeqn  \af_* ([1 \otimes z]) = n [z \otimes 1] + [1 \otimes z].
\]
To calculate $\af_* (\bt)$, we shift to the topological K-theory
$K^* (S^1 \times S^1)$, and use its ring structure.
Note that $[z \otimes 1]^2 = 0$, because all elements of $K^1 (S^1)$
have square zero.
Therefore
\[
\af_* (\bt) =
 {\ts{ \left( \rsz{ h^{-1} } \right)^* }}
                ( [z \otimes 1] \cdot [1 \otimes z] )
 = [z \otimes 1] \cdot (n [z \otimes 1] + [1 \otimes z]) = \bt.
\]
It follows that
$\af_* \colon K_0 ( C ( S^1 \times S^1 ) ) \to K_0 ( C ( S^1 \times S^1 ) )$
is the identity, and that
$\af_* \colon K_1 ( C ( S^1 \times S^1 ) ) \to K_1 ( C ( S^1 \times S^1 ) )$
is given by the matrix
{\scriptsize{$
 \left( \begin{array}{cc} 1 & n \\ 0 & 1 \end{array} \right) $}}.

We now know that the upper right horizontal map in the
Pimsner-Voiculescu sequence is zero, and that the lower left
horizontal map is
{\scriptsize{$
 \left( \begin{array}{cc} 0 & - n \\ 0 & 0 \end{array} \right) $}}.
Therefore
\[
K_1 (A) \cong \Z^3 \oplus \Z / n \Z \andeqn K_0 (A) \cong \Z^3.
\]
We identify $[1]$ and $\bt$ with their images
in $K_0 (A)$, and choose any $\et_0 \in K_0 (A)$ such that
$\partial (\et_0) = (-1, 0) \in K_1 ( C ( S^1 \times S^1 ))$.
(Note that $(-1, 0)$ generates the kernel of $\id - \af_*^{-1}$ on
$K_1 ( C ( S^1 \times S^1 ) )$.)
Then $K_0 (A) = \Z \cdot [1] \oplus \Z \cdot \bt \oplus \Z \cdot \et_0$.
(We will replace $\et_0$ by a more carefully chosen generator later.)

Now let $\ta$ be an arbitrary normalized trace on $A$.
Then $\ta$ is induced by
an $h$-invariant measure $\mu$ on $S^1 \times S^1$.
We compute $\ta_* \colon K_0 (A) \to \R$.
Trivially $\ta (1) = 1$,
and $\ta_* (\bt) = 0$ because $\bt \in K^0 (S^1 \times S^1)$
is represented as the difference of two vector bundles of the same rank
(namely $1$).
To calculate $\ta_* (\et_0)$,
we combine Definition~VI.8 and Theorems~V.12 and~VI.11 of \cite{Ex}
to get (notation explained afterwards)
\[
\exp (2 \pi i \ta_* (\et_0) ) = R_{\af}^{\mu} ([z^{-1} \otimes 1]).
\]
Here $[z^{-1} \otimes 1]$ now represents the homotopy class of the
function $(\zt_1, \zt_2) \mapsto \zt_1^{-1}$
(an element of $[S^1 \times S^1, \, S^1]$).
Following Definitions~VI.3 and~VI.5 of \cite{Ex},
$R_{\af}^{\mu} ([v])$ is computed by finding a
\ct\  function  $f \colon S^1 \times S^1 \to \R$ such that
\[
{\ts{ v \left( \rsz{ h^{-1} (x) } \right)^* v (x) }} = e^{i f (x)}
\]
for all $x \in S^1 \times S^1$, and setting
\[
R_{\af}^{\mu} ([v]) = \exp \left( i \int_X f \, d \mu \right).
\]
(By comparing Definition~VI.2 with the proof of Proposition~VI.10
in \cite{Ex}, one sees that the automorphism $C ( S^1 \times S^1 )$
given by $h$ really is $\af (f) = f \circ h^{-1}$.)
With $v = z \otimes 1$, one checks we may choose the function
$f (x) = 2 \pi \te$ for all $x$,
whence $\exp (2 \pi i \ta_* (\et_0) ) = \exp (2 \pi i \te )$.
Therefore there is $k \in \Z$ such that $\ta_* (\et_0) = \te + k$.

A priori $k$ depends on $\ta$.
However, the space of normalized traces is connected, and
$\ta \mapsto \ta_* (\et_0)$ is \ct, so in fact $k$ is independent
of $\ta$.
Replacing $\et_0$ by $\et_0 - k [1]$, we may therefore assume
that $\ta_* (\et_0) = \te$ for all traces $\ta$.
It follows that $\ta_* \colon \Z^3 \to \R$ is given, for every $\ta$, by
$\ta_* (m_1, m_2, m_3) = m_1 + m_3 \te$.

Every normalized trace $\ta$ on $A$ must restrict to the unique
trace on the image $\ph (A_{\te}) \subset A$ of the irrational rotation
algebra $A_{\te}$.
Therefore, with the \pj\  $p \in A_{\te}$ being as at the beginning,
we have $\ta_* ([ \ph (p)] - \et_0) = 0$.
Consequently there is $l \in \Z$ such that $\et_0 = [ \ph (p)] + l \bt$.
We then also have
\[
K_0 (A) = \Z \cdot [1] \oplus \Z \cdot \bt \oplus \Z \cdot [ \ph (p)],
\]
and the formula for every $\ta_*$ with respect to the new
identification of $K_0 (A)$ with $\Z^3$ is still
$\ta_* (m_1, m_2, m_3) = m_1 + m_3 \te$.
The identification of $K_0 (A)_+$ with
\[
\{ (m_1, m_2, m_3) \in \Z^3 \colon m_1 + m_3 \te > 0 \,\, {\text{or}} \,\,
  m_1 = m_2 = m_3 = 0 \}
\]
(the main result of \cite{Kd})
is now immediate from Theorem~\ref{Q5}~(2).
\end{exa}

\end{document}